\input amstex

\documentstyle{amsppt}
\nologo
\NoBlackBoxes

\magnification\magstep1
\pagewidth{32.622pc}
\pageheight{44.8pc}

\def\Ad{\operatorname{Ad}}
\def\Hom{\operatorname{Hom}}

\def\Aut{\operatorname{Aut}}
\def\Int{\operatorname{Int}}
\def\Ind{\operatorname{Ind}}
\def\ind{\operatorname{ind}}

\def\pro{\operatorname{pro}}
\def\Ext{\operatorname{Ext}}

\def\GL{\operatorname{GL}}

\def\ad{\operatorname{ad}}

\def\coker{\operatorname{coker}}

\def\For{\operatorname{For}}

\def\deg{\operatorname{deg}}

\def\A{\Cal A}
\def\B{\Cal B}

\def\D{\Cal D}
\def\E{\Cal E}
\def\U{\Cal U}
\def\M{\Cal M}
\def\Sl{\Cal S}

\def\N{\Cal N}
\def\I{\Cal I}

\def\d{\frak d}
\def\g{\frak g}
\def\k{\frak k}

\def\t{\frak t}
\def\p{\frak p}

\def\Vd{V^\cdot}

\def\Homd{\Hom^\cdot}

\def\ootimes{\overline{\otimes}}

\def\eps{\varepsilon}

\def\pf{\demo{Proof}}

\topmatter

\title
Equivariant analogues of Zuckerman functors
\endtitle
\author 
Pavle Pand\v zi\'c
\endauthor
\affil 
Department of Mathematics,
University of Zagreb, Croatia  
\endaffil
\address 
Department of Mathematics,
University of Zagreb,
PP 335, 10002 Zagreb, Croatia
\endaddress
\email 
pandzic\@math.hr 
\endemail
\abstract We review the Beilinson-Ginzburg construction of equivariant
derived categories of Harish-Chandra modules, and introduce analogues
of Zuckerman functors in this setting. They are given by
an explicit formula, which works equally well in the case of modules with
a given infinitesimal character. This is important if one wants to
apply Beilinson-Bernstein localization, as is explained in \cite{MP1}.
We also show how to recover the usual Zuckerman functors from the
equivariant ones by passing to cohomology.
\endabstract
\endtopmatter

\document

\head Introduction
\endhead

One of the basic problems in representation theory is to understand
linear representations of a real semisimple connected Lie group $G_0$
with finite center. One can instead study related problems of
understanding $\frak g$-modules, where $\frak g$ is the complexified
Lie algebra of $G_0$, or $(\frak g,K)$-modules, where $K$ is the
complexification of a maximal compact subgroup of $G_0$.

To be more specific, let $\Cal M(\g)$ be the abelian
category of modules over $\g$, or equivalently over the enveloping 
algebra $\Cal U(\frak g)$ of
$\frak g$. Let $K$ be an algebraic group acting on $\frak g$ via inner
automorphisms, such that the differential of this action identifies
the Lie algebra $\frak k$ of $K$ with a subalgebra of $\frak g$. Then
one can define an abelian category $\Cal M(\g,K)$ of {\it
Harish-Chandra modules} for the pair $(\frak g,K)$. This category
consists of $\g$-modules with an algebraic action of $K$
satisfying the usual compatibility conditions; see \S 1.
Clearly, there is a natural forgetful functor from 
$\Cal M(\g,K)$ into $\Cal M(\g)$. 

In the late seventies, G.~Zuckerman proposed a
homological construction of Harish-Chandra modules based on the
following observation. Let $T$ be a closed subgroup of $K$. Then there
is a natural forgetful functor from $\Cal M(\g,K)$
into $\Cal M(\g,T)$. This functor has a right adjoint
functor $\Gamma_{K,T}$, called the {\it Zuckerman functor}. Roughly 
speaking, 
$\Gamma_{K,T}$ attaches to a Harish-Chandra module in 
$\Cal M(\g,T)$ its largest Harish-Chandra submodule for $(\frak
g,K)$. Zuckerman functors are left exact and usually zero on
``interesting'' modules. Therefore, one uses the machinery of
homological algebra and considers their right derived functors.
Using them, one can construct ``a lot of'' $(\g,K)$-modules. For example,
if $\g$ and $K$ are related to $G_0$ as above, then
all irreducible $(\g,K)$-modules are submodules of modules obtained
in this way for certain specific choices of $(\g,T)$-modules for
various suitable (reductive) $T$. See for example \cite{Vo} or \cite{KV}.
We review the Zuckerman construction briefly in \S 1; a more detailed
review is contained in \cite{MP1}, \S 1. 

As usual, if one wants to study homological algebra of an abelian
category $\Cal A$, the natural setting is the derived category $D(\Cal
A)$ of $\Cal A$. In our case, one can consider the bounded derived
categories $D^b(\Cal M(\g,T))$ and 
$D^b(\Cal M(\g,K))$ of complexes of Harish-Chandra modules, the forgetful
functor from $D^b(\Cal M(\g,K))$ into 
$D^b(\Cal M(\g,T))$ and its right adjoint, the {\it derived Zuckerman
functor} $R\Gamma_{K,T}$ from $D^b(\Cal M(\g,T))$ into
$D^b(\Cal M(\g,K))$. One can recover the ``classical'' derived Zuckerman 
functors by applying $R\Gamma_{K,T}$ to a single module and then taking
cohomology.

On the other hand, a powerful algebro-geometric method of studying
$\frak g$-modules is the localization theory, due to A.~Beilinson and 
J.~Bernstein. Let $\theta$ be a maximal ideal in the
center $\Cal Z(\frak g)$ of the enveloping algebra $\Cal U(\frak g)$
and denote by $\Cal U_\theta$ the quotient of $\Cal U(\frak g)$ by the
two-sided ideal generated by $\theta$. As before, we can define the
abelian categories $\Cal M(\Cal U_\theta)$ and $\Cal M(\Cal
U_\theta,K)$ of $\g$-modules, resp.~Harish-Chandra
modules, with a given infinitesimal character determined by
$\theta$. By the Harish-Chandra homomorphism, $\theta$ corresponds to a
Weyl group orbit in the linear dual $\frak h^*$ of a Cartan subalgebra
$\frak h$ of $\frak g$. To each $\lambda$ in $\frak h^*$, Beilinson
and Bernstein attach a {\it sheaf of twisted differential operators}
$\Cal D_\lambda$ on the flag variety $X$ of $\frak g$.  They prove
that the global sections $\Gamma(X,\Cal D_\lambda)$ of $\Cal
D_\lambda$ are equal to $\Cal U_\theta$, and that the higher cohomology
vanishes.  They define the localization functor from the category $\Cal
M(\Cal U_\theta)$ into the category of quasicoherent $\Cal
D_\lambda$-modules $\Cal M(\Cal D_\lambda)$ on $X$ by
$\Delta_\lambda(V) =\Cal D_\lambda \otimes_{\Cal U_\theta} V$. In the
case of a regular infinitesimal character, and $\lambda$ which
satisfies a geometric positivity condition, the localization functor
$\Delta_\lambda$ is an equivalence of the category $\Cal M(\Cal
U_\theta)$ with $\Cal M(\Cal D_\lambda)$. Its quasi-inverse is the
functor of global sections $\Gamma(X, - )$. Under this equivalence,
the category $\Cal M(\Cal U_\theta,K)$ of Harish-Chandra modules
corresponds to the category $\Cal M(\Cal D_\lambda,K)$ of {\it
Harish-Chandra sheaves}.  Without the positivity condition, the
localization functor induces an equivalence of corresponding derived
categories $D^b(\Cal M(\Cal U_\theta))$ and $D^b(\Cal M(\Cal
D_\lambda))$. This enables one to apply the machinery of sheaf
cohomology theory to the study of representations. Many results have
been obtained using this method, the most famous one being the proof
of Kazhdan-Lusztig conjectures.

It is a natural problem to try to connect these two
constructions and ``localize'' the Zuckerman functor
construction. The solution is not obvious, as the homological properties 
of the categories $\Cal M(\g,K)$ and $\Cal M(\Cal D_\lambda)$
are quite different. Still, Bernstein managed to find a partial solution 
in the case when $\lambda$ satisfies the
geometric positivity condition so that the localization functor is
exact. Also, Hecht, Mili\v ci\'c, Schmid and Wolf in their Duality
Theorem calculated the cohomology of ``standard'' Harish-Chandra
sheaves on the flag variety in terms of ``standard'' Zuckerman
modules, which can be interpreted as a particular instance of this
problem.

In the meantime, the notion of equivariant derived category of Harish-Chandra modules
(and D-modules) has been introduced by Beilinson and Ginzburg (\cite{BB}, \cite{G}; 
see also \cite{BL2} and \cite{MP1}, \S 2). An analogous construction for
equivariant sheaves has been done by Bernstein and Lunts. This category can
be used as a bridge between the above two theories. Namely, it allows localization
in a natural way, and on the other hand it carries the information about 
homological algebra of  Harish-Chandra modules as well. The reason for the
better behavior of this category with respect to localization is the fact
that its
definition is based on the notion of ``weak Harish-Chandra modules.''
These are defined in the same way as Harish-Chandra modules, but the
compatibility conditions are relaxed: the two actions of $\k$ are no longer
required to agree. In particular, free $\U(\g)$-modules are weak Harish-Chandra
modules for the pair $(\g,K)$ (but they are not Harish-Chandra modules). 
Similarly, free $\U_\theta$-modules are 
weak $(\U_\theta,K)$-modules, and consequently it is easy to study 
localization theory for such modules. In \S 2 we review the 
Beilinson-Ginzburg construction, including some details that are usually 
not mentioned.

A technical complication is that the equivariant derived category of
Harish-Chandra modules $D^b(\g,K)$ is a triangulated
category, but not a priori a derived category of any abelian
category. Therefore, its use requires some modifications of standard
tools used in homological algebra. In particular, 
a generalization of the notion of derived functors to the case of an
arbitrary
triangulated category and its localization due to P.~Deligne, 
\cite{De} is used. In this setting the place of projective and 
injective resolutions is taken by the notions of
$\Sl$-projective and $\Sl$-injective objects in the
triangulated category in question; here $\Sl$ is the localizing class defining
the localization. These objects were first defined by Verdier 
in \cite{Ve}; he calls them ``free on the left'' (respectively right).
In the setting of homotopic categories, 
they were further studied and applied by Spaltenstein in \cite{Sp} and
Bernstein and Lunts in \cite{BL1} and \cite{BL2} under the name of 
K-projectives and K-injectives.
One can find an account of this material in the Appendix of \cite{P2},
together with a summary of some mostly well-known properties of adjoint
functors, with emphasis on their use in homological algebra. 

\S 3 is devoted to the construction of the equivariant analogues 
$\Gamma_{K,T}^{eq}$  of the derived Zuckerman functors. $\Gamma_{K,T}^{eq}$ 
is defined as the right adjoint of the forgetful 
functor from the equivariant derived category $D(\frak g,K)$ 
into $D(\frak g,T)$. 
Analogously, one can define a variant of this functor
$\Gamma_{K,T}^{eq} : D(\Cal U_\theta,T) @>>> D(\Cal U_\theta,K)$.
It is shown how both of them can be described by the same explicit
formula, involving the standard complex $N(\k)$ of $\k$. The main 
technical result here is that $N(\k)$
is a K-projective equivariant $(\k,T)$-complex; see \S 3.2.
Furthermore, it is proved that in the case of 
$\Gamma_{K,T}^{eq} : D(\g,T) @>>> D(\g,K)$, when applied to a single 
Harish-Chandra module $V$, this formula gives a complex with cohomology
modules equal to the ``classical'' derived Zuckerman functors of $V$. 
In fact, in this 
case the formula reduces to the Duflo-Vergne formula 1.3.3; see
\cite{DV} and \cite{MP1}, 1.6. All of the above is also true on the
level of the bounded equivariant derived categories.

It is explained in \cite{MP1} how these results can be used to localize
the Zuckerman construction.

Finally, let me mention that, as indicated in \cite{MP1}, one can use
localization of the Zuckerman construction to obtain a formula for
the cohomology of standard Harish-Chandra sheaves in terms of the derived
Zuckerman functors. As a special case one can recover the Duality theorem
of Hecht, Mili\v ci\'c, Schmid and Wolf from \cite{HMSW}. The details
will appear in \cite{MP2}.

Most of this paper comes out of my 1995 University of Utah Ph.D. thesis,
which was written under guidance of my advisor Dragan Mili\v ci\'c.
I am very much indebted to him for his generous help and many suggestions
and ideas he supplied me with. I would also like to thank David Vogan for several
useful conversations.

\head 1. Harish-Chandra modules and Zuckerman functors 
\endhead

\subhead 1.1. Harish-Chandra pairs and modules \endsubhead
Let $K$ be a complex algebraic group and $\frak k$ its Lie algebra. 
We recall that an algebraic representation of $K$ is a vector space $V$
with a $K$-action, which is a union of $K$-invariant finite-dimensional
subspaces $V_i$, so that the action on each $V_i$ is given by an algebraic 
group morphism $K@>>> \GL(V_i)$.

Let $\Cal A$ be an associative algebra with identity, with an algebraic action 
of $K$ 
$$
\phi:K@>>>\Aut(\Cal A),
$$
and a Lie algebra homomorphism $\psi$ from $\frak k$ into $\Cal A$ such that
 the following conditions hold:
\roster
\item"(A1)" $\psi$ is $K$-equivariant, i.e.,
$$
\psi\circ\Ad k  = \phi(k) \circ\psi \text{ for all } k \in K;
$$
\item"(A2)" The differential of the action $\phi$ is the same as the action
of $\frak k$ on $\Cal A$ given by inner derivations $\psi(\xi)$, $\xi \in \frak k$, i.e., 
$$
L(\phi)(\xi)(a)=[\psi(\xi),a]
$$
for any $\xi \in\frak k$ and $a\in\Cal A$.
\endroster
These conditions are related, but in a general situation none of
them implies the other.

A pair $(\A,K)$ as above is called a Harish-Chandra pair.
The main examples come from a similar notion of a ``classical" Harish-Chandra
pair, when $\A$ is replaced by a Lie algebra $\g$. Let $\g$ be a complex
Lie algebra, and $K$ a (usually reductive) algebraic group over $\Bbb C$, acting 
algebraically on $\g$ via a morphism $\phi:K@>>>\Int(\g)$ of algebraic groups.
Assume further there is a $K$-equivariant embedding $\psi$ of the Lie algebra 
$\k$ of $K$ into $\g$, such that the differential of $\phi$ is equal to 
$\ad\circ\psi$. Then $(\g,K)$ is called a (classical) Harish-Chandra pair.

In this situation, $(\U(\g),K)$ is a Harish-Chandra pair as above: $\phi$
defines an action of $K$ on $\U(\g)$ in an obvious way, and $\psi$ defines
an embedding of $\k$ into $\U(\g)$; the compatibility conditions (A1) and (A2)
clearly hold. An important special case of this is when $G$ is a real semisimple 
Lie group, $\g$ its complexified Lie algebra, and $K$ the complexification of 
the maximal compact subgroup of $G$. The map $\phi$ is just the adjoint action 
and $\psi$ is the inclusion of $\k$ into $\U(\g)$.

Another set of examples is obtained from a complex Lie group $G$; let $\g$ be 
the Lie algebra of $G$, and $K$ be any algebraic
subgroup of $G$, like a Borel subgroup or its unipotent radical; $K$ can also
be $G$ itself. The maps $\phi$ and $\psi$ are, as before,
the adjoint action and the inclusion.

In the above cases, when $\g$ is semisimple, one can also 
take $\A$ to be the quotient $\U_\theta$ of $\U(\g)$ corresponding 
to an infinitesimal character given by the Weyl group orbit $\theta$ of an
element of the dual of a Cartan subalgebra of $\g$. This is particularly important
for us, as it is the setting of the Beilinson-Bernstein localization theory.

Finally, we will occasionally us another simple example: $\A=\Bbb C$, and $K$ 
is an arbitrary algebraic group;
$\phi$ is the trivial action on $\Bbb C$ and $\psi$ is 0.

A vector space $V$ is called a {\it weak $(\Cal A,K)$-module} or a weak 
Harish-Chandra module for the pair $(\A,K)$ if:
\roster
\item"(HC1)" $V$ is an $\Cal A$-module with an action $\pi$;
\item"(HC2)" $V$ is an algebraic $K$-module with an action $\nu$;
\item"(HC3)" for any $a \in \Cal  A$ and $k \in K$ we have
$$
\pi(\phi(k)a)=\nu(k)\pi(a)\nu(k)^{-1};
$$
i.e., the $\Cal A$-action map $\Cal A \otimes V @>>> V$ is $K$-equivariant.
\endroster
In future we will denote the representation $\pi\circ\psi$ of $\frak k$
on $V$ simply by $\pi$. This should not create confusion; in fact in
the most important examples $\psi$ is injective and $\frak k$ can be
identified with a Lie subalgebra of $\Cal A$.

The action $\nu$ of $K$ differentiates to an action of $\frak k$ which we 
denote also by $\nu$. We put $\omega(\xi)=\nu(\xi)-\pi(\xi)$ for
$\xi\in\frak k$. The following lemma is \cite{MP1, 2.1}.
\proclaim
{1.1.1. Lemma} Let $V$ be a weak $(\Cal A,K)$-module. Then 
\roster
\item"(i)" $\omega$ is a representation of $\frak k$ on $V$;
\item"(ii)" $\omega$ is $K$-equivariant, i.e.,
$$
\omega(\Ad(k)\xi) = \nu(k) \omega(\xi) \nu(k)^{-1}
$$
for $\xi \in \frak k$ and $k \in K$.
\item"(iii)" 
$$
[\omega(\xi),\pi(a)]=0,
$$
for $a \in \Cal A$ and $\xi \in \frak k$.
\qed\endroster
\endproclaim

We say that a weak $(\Cal A,K)$-module $V$ is an $(\Cal A,K)$-module if
$\omega=0$, i.e., if in addition to (HC1)-(HC3) it satisfies:
\roster
\item"(HC4)" $\pi(\xi) = \nu(\xi)$ for all $\xi\in\frak k$.
\endroster
A morphism $\alpha : V @>>> W$ of two weak $(\Cal A,K)$-modules is a linear 
map which is a morphism for both $\Cal A$- and $K$-module structures.
We denote by $\Cal M(\Cal A,K)_w$ the category of all weak
$(\Cal A,K)$-modules, and by $\Cal M(\Cal A,K)$ its full subcategory of
$(\Cal A,K)$-modules. Clearly, these are abelian categories. 

The inclusion functor from $\M(\A,K)$ into $\M(\A,K)_w$ has both adjoints.
They are defined as follows.
Let $U \in \Cal M(\Cal A,K)_w$. Let
$$
U^{\frak k} =\{ u \in U \mid \omega(\xi) u = 0, \ \xi \in \frak k\}
$$
be the $\k$-invariants in $U$ with respect to $\omega$ and let
$$
U_{\frak k}=U/\omega(\k)U
$$
be the $\k$-coinvariants of $U$ with respect to $\omega$. Then 1.1.1 implies 
that the actions of $\A$ and $K$ on $U$ induce actions on $U^\k$ and $U_\k$. 
It is now clear that $U^\k$ is the largest $(\A,K)$-submodule of $U$, and that
$U_\k$ is the largest $(\A,K)$-quotient of $U$. This immediately implies
the following proposition.

\proclaim
{1.1.2. Proposition} For any $(\A,K)$-module $V$
$$
\Hom_{(\A,K)}(V,U)=\Hom_{(\A,K)}(V,U^\k) \quad \text{and}
$$
$$
\Hom_{(\A,K)}(U,V)=\Hom_{(\A,K)}(U_\k,V).
$$
In other words, the functors $U \longmapsto U^{\frak k}$ and
$U \longmapsto U_{\frak k}$ are right respectively left adjoint to
the forgetful functor from $\Cal M(\Cal A,K)$ into $\Cal M(A,K)_w$.
\qed\endproclaim

\subhead 1.2. Change of algebras \endsubhead
Let $(\Cal A,K)$ and $(\Cal B,K)$ be two Harish-Chandra pairs. Denote by 
$\phi_{\Cal A}$, $\psi_{\Cal A}$ and $\phi_{\Cal B}$, $\psi_{\Cal B}$ the
corresponding maps. Assume that there is a morphism of algebras 
$\gamma : \Cal B @>>> \Cal A$ such that
$$
\gamma \circ \phi_{\Cal A}(k) = \phi_{\Cal B}(k) \circ \gamma, \quad k \in K,
$$ 
i.e., $\gamma$ is $K$-equivariant.  If $V$ is a weak $(\Cal A,K)$-module, we 
can view it as a weak $(\Cal B,K)$-module, where the
action of $\Cal B$ is given by $\tilde\pi(b) = \pi(\gamma(b))$ for $b
\in \Cal B$. The equivariance condition is preserved because of the
equivariance of $\gamma$. This defines an exact functor $\For =
\For_{\Cal B,\Cal A}$ from the category $\Cal M(\Cal A,K)_w$ into the
category $\Cal M(\Cal B,K)_w$. This forgetful functor has both adjoints.
We will need the left adjoint, which we now describe.

Let $V$ be a weak $(\B,K)$-module. We can consider $\A$ as a right $\B$-module
via the right multiplication composed with $\gamma$. Then we can form the
vector space
$$
\ind_{\A,\B}(V)=\A\otimes_\B V.
$$
One easily checks that this space becomes a weak $(\A,K)$-module if we
let $\A$ act by left multiplication in the first factor, and $K$ by the
tensor product of the actions $\phi_\A$ on $\A$ and $\nu_V$ on $V$.
Furthermore, $\ind_{\A,\B}$  is a functor from $\M(\B,K)_w$ into $\M(\A,K)_w$.
If we forget the $K$-action, then it is a well-known fact that this functor
is left adjoint to the forgetful functor from $\A$-modules to $\B$-modules.
However, the standard adjunction morphisms 
$v\mapsto 1\otimes v,\,\, v\in V$ for a $\B$-module $V$, and 
$a\otimes v\mapsto \pi_V(a)v,\,\, a\in\A,\,v\in V$ for an 
$\A$-module $V$, are easily checked to
be $K$-equivariant in our present situation. This implies the following
proposition.

\proclaim{1.2.1. Proposition} Let $\gamma:\B@>>>\A$ be a K-equivariant 
morphism of algebras. Let $\For$ be the corresponding forgetful functor 
from $\M(\A,K)_w$ into $\M(\B,K)_w$, and let the functor $\ind_{\A,\B}$ 
be defined as above. Then $\ind_{\A,\B}$ is left adjoint to $\For$.
\qed\endproclaim

In case $K$ is reductive, we can use this result to construct projectives in 
the categories $\M(\A,K)_w$ and $\M(\A,K)$. Namely, clearly 
$\M(K)=\M(\Bbb C,K)_w$ and this category is semisimple, hence all its objects
are projective. Furthermore, the inclusion $\Bbb C@>>>\A$ given by 
$\lambda\mapsto\lambda 1$ is clearly $K$-equivariant. 
Since $\For_{\Bbb C,\A}$ is exact, its left
adjoint $\ind_{\A,\Bbb C}$ preserves projectives. So for any $V\in\M(\A,K)_w$,
$\ind_{\A,\Bbb C}(V)$ is a projective weak $(\A,K)$-module. On the other hand,
the adjunction morphism $\Psi_V:\ind_{\A,\Bbb C}(V)@>>>V$ is surjective, since
$\Psi_V(1\otimes v)=v$ for any $v\in V$. So $V$ is a quotient of a
projective module.

Using this together with 1.1.2, we can similarly prove that $\M(\A,K)$ also
has enough projectives. Hence we have proved

\proclaim{1.2.2. Proposition} Let $K$ be a reductive algebraic group. Then the 
categories $\M(\A,K)_w$ and $\M(\A,K)$ have enough projectives.
\qed\endproclaim

The other adjoint is constructed similarly. It is given by
$$
\pro_{\A,\B}(V)=\Hom_\B(\A,V)^{alg}
$$
for a weak $(\B,K)$-module $V$, with suitably defined actions. Here $W^{alg}$
denotes the $K$-algebraic part of a $K$-module $W$ which is not necessarily
algebraic. We omit the details, since we do not need the functor $\pro_{\A,\B}$
 in this paper.

Let us note that if we in addition assume that $\psi_\A=\gamma\circ\psi_\B$,
then clearly $\For_{\B,\A}$ sends $\M(\A,K)$ into $\M(\B,K)$. Moreover, an
easy calculation shows that the $\omega$-action on $\ind_{\A,\B}V$ is equal to 
the tensor product of the trivial action on $\A$ with $\omega_V$.
It follows that $\ind_{\A,\B}$ sends $(\B,K)$-modules into $(\A,K)$-modules.
This implies the following variant of 1.2.1.

\proclaim{1.2.3. Proposition} Let $\gamma:\B@>>>\A$ be a K-equivariant 
morphism of algebras such that $\psi_\A=\gamma\circ\psi_\B$. Then the
functor $\ind_{\A,\B}:\M(\B,K)@>>>\M(\A,K)$ is left adjoint to the forgetful
functor from $\M(\A,K)$ into $\M(\B,K)$.
\qed\endproclaim

An analogous result holds for the functor $\pro_{\A,\B}$.

\subhead 1.3. Change of groups and Zuckerman functors \endsubhead
Let $(\Cal A,K)$ be a Harish-Chandra pair, with $\phi_K : K @>>> \Aut(\Cal A)$ 
and $\psi_K : \frak k @>>> \Cal A$ the corresponding morphisms.

Let $T$ be another algebraic group and $\gamma : T @>>> K$ a morphism of 
algebraic groups. The main example is the inclusion of a closed subgroup. 
We denote the differential of $\gamma$ again by $\gamma$; 
in applications it will always be injective.
It is easily checked that $(\Cal A,T)$ becomes a Harish-Chandra pair if
we define $\phi_T = \phi_K \circ \gamma$ and $\psi_T = \psi_K \circ \gamma$.
Moreover, composing the $K$-action with $\gamma$ to get a $T$-action clearly
defines a forgetful functor from $\Cal M(\Cal A,K)_w$ into 
$\Cal M(\Cal A,T)_w$, and also from $\Cal M(\Cal A,K)$ into $\Cal M(\Cal A,T)$.
 We want to describe right adjoints of these forgetful functors.

First we treat the case $T=\{1\}$. Then, as was proved in \cite{MP1, 2.2}
(see also \cite{MP1, 1.2}), the right adjoint of 
$\For:\Cal M(\Cal A,K)_w @>>> \Cal M(\Cal A)$ is $\Ind_w$, given by
$$
\Ind_w(V)=R(K)\otimes V = R(K,V).
$$
Here $V$ is an $\A$-module, $R(K)$ denotes the space of regular functions on 
$K$, and $R(K,V)$ denotes the space of $V$-valued regular functions on $K$. 
Recall that $K$ acts 
on $\Ind_w(V)$ by the right regular representation, and $a\in \A$ acts on 
$F\in\Ind_w(V)$ by
$$
(\pi(a)F)(k)=\pi_V(\phi_K(k)a)(F(k)),\quad k\in K.
$$
The adjunction morphisms are given as follows. For $V\in\Cal M(\Cal A,K)_w$,
$\Phi_V:V@>>>\Ind_w(\For V)$ is the matrix coefficient map, given by
$$
\Phi_V(v)(k)=\nu_V(k)v,\quad v\in V,\,k\in K.
$$
For $W\in\M(\A)$, $\Psi_W:\For(\Ind_w(W))@>>>W$ is the evaluation at $1$.

Composing $\Ind_w$ with the functor of taking $\k$-invariants with
respect to $\omega$ from 1.1.2, we get the functor $\Gamma_K$ which is
right adjoint to $\For:\Cal M(\Cal A,K)@>>> \Cal M(\Cal A)$. One can
use the functors $\Ind_w$ and $\Gamma_K$ to show the existence of enough 
injectives in $\Cal M(\Cal A,K)_w$ and $\Cal M(\Cal A,K)$; see \cite{MP1, 
1.3, 2.3}. Note that this is analogous to 1.2.2.

Now we go back to the `relative' case, when the group $T$ is nontrivial.
We first want to construct the right adjoint of
$\For:\M(\A,K)_w @>>> \M(\A,T)_w$. Let $V$ be a weak $(\A,T)$-module.
We can consider $V$ as an $\A$-module and construct the weak $(\A,K)$-module 
$\Ind_w(V)=R(K)\otimes V$ as above. We define an action $\lambda_T$ of $T$ 
on $\Ind_w(V)$ as the tensor product of the 
left regular representation of $T$ on $R(K)$ (via $\gamma$) with the given 
action of $T$ on $V$. Then it is easy to see that $\lambda_T$ commutes with
the actions $\pi$ and $\nu$ of $\A$ and $K$.

It follows that the space $I_w(V)$ of $\lambda_T$-invariants in $\Ind_w(V)$
is a weak $(\A,K)$-submodule of $\Ind_w(V)$. It is now easy to see that
$I_w:\M(\A,T)_w @>>> \M(\A,K)_w$ is the right adjoint of the forgetful
functor; the adjunction morphisms are obtained by restriction from the
adjunction morphisms for $\Ind_w$.

Taking invariants with respect to $\omega$, we get
the functor $V\mapsto I_w(V)^\k$, which is right adjoint to the forgetful 
functor $\For:\M(\A,K)@>>>\M(\A,T)_w$.
Since $\M(\A,T)$ is a full subcategory of $\M(\A,T)_w$, by restricting the
domain of this functor we get the right adjoint of the forgetful functor
$\For:\M(\A,K)@>>>\M(\A,T)$; see \cite{P2}, A2.1.1. This functor is denoted
by $\Gamma_{K,T}$, and called the {\it Zuckerman functor}.

Another construction of $\Gamma_{K,T}$ was explained in \cite{MP1, 
\S 1}; instead of $\k$-invariants with respect to the action $\omega$, we can
take $\k$-invariants with respect to the action $\lambda_\k$, which is the
tensor product of the left regular action of $\k$ on $R(K)$ with the given 
action on $V$. 
These two constructions agree because of the following formula which was proved
below 2.4 in \cite{MP1}, and which implies that taking $\k$-invariants
with respect to $\omega$ or $\lambda_\k$ gives the same result.

\proclaim{1.3.1. Lemma} For any $\xi\in\k$, $F\in\Ind_w(V)$ and $k\in K$,
$$
(\omega(\xi)F)(k)=-(\lambda_\k(\Ad(k)\xi)F)(k). \qed
$$
\endproclaim
So we have

\proclaim{1.3.2. Theorem} The forgetful functor from $\M(\A,K)$ to $\M(\A,T)$
has a right adjoint, the Zuckerman functor $\Gamma_{K,T}$. For an 
$(\A,T)$-module $V$, 
$$
\Gamma_{K,T}(V)=\Hom_{\k,T}(\Bbb C,R(K,V)),
$$ 
where $\Bbb C$ is the trivial $(\k,T)$-module, $R(K,V)$ is a
$(\k,T)$-module with respect to the action
$\lambda=(\lambda_\k,\lambda_T)$, and the $(\A,K)$-action on
$\Gamma_{K,T}(V)$ comes from the $(\A,K)$-action on $R(K,V)$.
\qed\endproclaim

Being a right adjoint, $\Gamma_{K,T}$ is left exact, and it has right derived 
functors (since there are enough injectives).
The following is a version of the Duflo-Vergne formula \cite{DV}, which 
is an explicit formula for the derived Zuckerman functors holding
under certain assumptions. It generalizes the formula for $\Gamma_{K,T}$ from
1.3.2. A proof can be found in \cite{MP1, 1.6}.

\proclaim{1.3.3. Theorem} Assume that $T$ is reductive and that $\A$ is a flat 
right $\U(\k)$-module for the right multiplication. Let $V$ be an 
$(\A,T)$-module. Then
$$
R^p \Gamma_{K,T}(V) = H^p(\k,T;\Ind_w(V))
$$
for $p \in \Bbb Z_+$, where $R(K,V)$ is a $(\k,T)$-module
with respect to the action $\lambda$. The $(\A,K)$-action on
$R^p \Gamma_{K,T}(V)$ comes from the $(\A,K)$-action on $R(K,V)$.
\qed\endproclaim

One can also consider the derived functor $R\Gamma_{K,T}$ between the
corresponding derived categories. However, 1.3.3 does not generalize readily
to a claim about this functor. The situation is better in the setting of
equivariant derived categories which we define in the following.

\head 2. Equivariant derived categories
\endhead
We are going to outline an algebraic construction of equivariant derived 
categories due to Beilinson and Ginzburg (see \cite{BB}, \cite{G} and \cite{BL2}). 
There is also a geometric construction due to Bernstein and Lunts \cite{BL1}.

\subhead 2.1. Equivariant complexes
\endsubhead
Let $(\A,K)$ be a Harish-Chandra pair. An {\it equivariant $(\A,K)$-complex}
is a complex $V=\Vd$ of weak $(\A,K)$-modules, together with a family of linear 
maps $i_{\xi}:V @>>> V$ of degree $-1$, depending linearly on
$\xi\in\k$, and satisfying the following conditions:
$$
i_{\xi}\pi(a)=\pi(a)i_{\xi},\quad \xi\in\k,\, a\in\A; \tag{EC1} 
$$
\vskip -1pc
$$
i_{\Ad(k)\xi}=\nu(k)i_{\xi}\nu(k)^{-1},\quad\xi\in\k,\, k\in K; \tag{EC2}
$$
$$
i_{\xi}i_{\eta}+i_{\eta}i_{\xi}=0,\quad \xi,\eta\in\k; \tag{EC3}
$$
$$
di_{\xi}+i_{\xi}d=\omega(\xi)\quad\xi\in\k. \tag{EC4}
$$
Here $\pi$ and $\nu$ are the actions of $\A$ and $K$ on $V$, $\omega=\nu-\pi$
is the action of $\k$ from \S 1.1, and $d$ is the differential of the complex $V$.

Note that for any $\xi\in\k$, $\omega(\xi)$ is a morphism of complexes of
$\A$-modules, which is homotopic to zero, $i_\xi$ being the
homotopy. It follows that $H^p(\omega(\xi))=0$ for any $p\in\Bbb Z$, so the
cohomology modules of $V$ are $(\A,K)$-modules.

A morphism between two equivariant $(\A,K)$-complexes is a morphism of 
complexes of weak $(\A,K)$-modules which commutes with the $i_\xi$'s.
We denote the category of equivariant $(\A,K)$-complexes by $C(\A,K)$. It
is clearly an abelian category.

\subhead 2.2. DG algebras and modules
\endsubhead
It will be convenient to reformulate the definition of equivariant complexes
using the notion of DG modules over DG algebras. For more about these see
\cite{Il}; see also \cite{BL1} and \cite{BL2}. The details needed here are 
written out in \cite{P1}.
The Lie algebra version of the theory is closely related with the theory of
Lie superalgebras, which can be looked up in \cite{Sch}. 

In the following, a complex will often be identified with the sum of its
components with the natural grading.
A differential graded (DG) Lie algebra is a complex of vector spaces $\d$ with
a graded bilinear operation $[-,-]$ called the supercommutator, satisfying the 
following conditions for any three homogeneous $X,Y,Z\in\d$:
$$
[X,Y]=-(-1)^{\deg X\deg Y}[Y,X]; \tag{DGLA1}
$$
$$
[X,[Y,Z]]=[[X,Y],Z]+(-1)^{\deg X\deg Y}[Y,[X,Z]]; \tag{DGLA2}
$$
$$
d[X,Y]=[dX,Y]+(-1)^{\deg X}[X,dY]. \tag{DGLA3}
$$
If we forget the differential and introduce the obvious $\Bbb Z_2$-grading,
we get a Lie superalgebra.

A differential graded (DG) module over $\d$ is a complex $V$ of vector spaces,
with a graded action of $\d$ on $V$ by linear endomorphisms, satisfying the 
properties
$$
[X,Y]v=XYv-(-1)^{\deg X\deg Y}YXv 
$$
and
$$
d_V(Xv)=(d_\d X)v+(-1)^{\deg X}Xd_Vv
$$
for any (homogeneous) $X,Y\in\d$ and $v\in V$. 

An example of a DG Lie algebra is a Lie algebra considered as a complex 
concentrated in degree zero. The example important for us is the DG Lie 
algebra $\underline{\k}$, corresponding to a Lie algebra $\k$, defined 
as follows. As a complex of vector spaces, $\underline{\k}$ is 
$$
\dots 0@>>>\k@>1>>\k@>>>0\dots,
$$ 
with $\k$ appearing in degrees $-1$ and $0$. The supercommutator is defined in 
the following way: on $\underline{\k}^0=\k$ it is just the commutator of $\k$, 
$\underline{\k}^{-1}$ supercommutes with itself, and for $X\in\underline{\k}^0$, 
$\xi\in \underline{\k}^{-1}$, $[X,\xi]$ in $\underline{\k}$ is the same as
$[X,\xi]$ in $\k$, understood as an element of $\underline{\k}^{-1}$.

It is straightforward to check that in this way $\underline{\k}$ becomes a
DG Lie algebra. Furthermore, $K$ acts on $\underline{\k}$ by adjoint action
on each component. This is clearly an algebraic action by chain maps.

The reason for our interest in this is the following.
Let $V$ be an equivariant $(\A,K)$-complex. We can define an action of
$\underline{\k}$ on $V$ as follows: $X\in\underline{\k}^0=\k$ acts by
$\omega(X)$ while $\xi\in \underline{\k}^{-1}$ acts by $i_\xi$. The fact that
this is a DG action is equivalent to (EC3), (EC4), and the fact that $\omega$ 
is a representation of $\k$. Furthermore, (EC1), (EC2), and the analogous
properties of $\omega$ from 1.1, say that this is a $K$-equivariant action
which commutes with the action of $\A$.

Now we want to pass to associative algebras. An associative DG algebra is
a graded algebra $\D$ with unit, which is also a complex of vector spaces,
such that $d(1)=0$ and for any homogeneous $x,y\in\D$,
$$
d(xy)=d(x)y+(-1)^{\deg x}xd(y).
$$
Any (associative) DG algebra is a DG Lie algebra, if we define
$$
[a,b]=ab-(-1)^{\deg a\deg b}ba
$$
for homogeneous $a$ and $b$ in $\D$.

A morphism of two DG algebras is a map that is both an algebra homomorphism 
and a morphism of complexes. A (left) DG module over a DG algebra $\D$ is a 
graded (left) $\D$-module $M$, which is also a complex of vector spaces, 
such that the grading of the module
is the same as the grading of the complex, and the following condition holds:
$$
d_M(xm)=(d_\D x)m+(-1)^{\deg x}x(d_M m)
$$
for any homogeneous $x\in\D$ and $m\in M$.

An example of a DG algebra is any associative algebra considered as a complex
concentrated in degree 0. The most important examples for our purposes are the
universal enveloping algebras of DG Lie algebras (concretely, of 
$\underline\k$). These are defined analogously to ordinary enveloping algebras,
and enveloping algebras of superalgebras (see \cite{Sch}). Starting with
a DG Lie algebra $\d$, we first define its tensor algebra $T(\d)$; as an 
algebra, it is the ususal tensor algebra. The grading is given by 
$$
\deg(X_1\otimes\dots\otimes X_k)=\deg X_1+\dots+\deg X_k
$$
for homogeneous $X_1,\dots ,X_k\in\d$, and the differential is given by
$$
d(X_1\otimes\dots\otimes X_k)=\sum_{i=1}^k (-1)^{\deg X_1+\dots\deg X_{i-1}}X_1\otimes\dots\otimes X_{i-1}\otimes d(X_i)\otimes X_{i+1}\otimes\dots\otimes X_k
$$
and $d(1)=0$. It is easy to check that in this way $T(\d)$ becomes a DG
algebra. Then we consider the two-sided ideal $\I$ of $T(\g)$
generated by the elements of the form
$$
X\otimes Y-(-1)^{\deg X \deg Y}Y\otimes X-[X,Y]
$$
for homogeneous $X,Y\in\d$. One checks that $\I$ is a DG ideal, i.e., a graded
ideal closed under the differential. It follows that the universal enveloping
algebra 
$$
\U(\d)=T(\d)/\I
$$
of $\d$ is a DG algebra. There is an obvious DG Lie algebra morphism 
$\phi:\d@>>>\U(\d)$, and the following universal property holds: for any DG 
algebra $\D$ and a DG Lie algebra morphism $f$ from $\d$ into $\D$ (considered 
as a DG Lie algebra as above), there is a unique DG algebra morphism $\tilde f$
 from $\U(\d)$ into $\D$ such that $f=\tilde f\circ\phi$.
In particular, $M$ is a DG module over $\d$ if and only if it is a DG
module over $\U(\d)$.

The proof is the same as for Lie algebras, as is the proof of 
the following version of the Poincar\'e-Birkhoff-Witt Theorem.

\proclaim{2.2.1. Proposition} Let $X_1,\dots,X_n,Y_1,\dots,Y_m$ be a basis 
of $\d$ 
consisting of homogeneous elements, with $X_k$'s even and $Y_k$'s odd. Then
$$
(X_1^{i_1}\dots X_n^{i_n}Y_1^{j_1}\dots Y_m^{j_m})_{i_k\in\Bbb Z_+,\,j_k\in\{0,1\}}
$$
is a basis of $\U(\d)$.
\qed\endproclaim

Using this and a straightforward calculation of the differential, we obtain
the following (known) result which goes back to Cartier; see \cite{BL2}, 1.9.7. 
Recall that the standard complex $N(\k)$ of a Lie algebra $\k$ is defined by
$$
N^p(\k)=\U(\k)\otimes \tsize{\bigwedge^{-p}}(\k),
$$
with the differential given by
$$
\multline
d(u\otimes\lambda_1\wedge\dots\wedge\lambda_p)=
\sum_{i=1}^p (-1)^{i-1}u\lambda_i\otimes\lambda_1\wedge\dots\wedge\hat{\lambda_i}\wedge\dots\wedge\lambda_p+\\
\sum_{j<i}(-1)^{j+i}u\otimes[\lambda_j,\lambda_i]\wedge\lambda_1\wedge\dots\wedge\hat{\lambda_j}\wedge\dots\wedge\hat{\lambda_i}\wedge\dots\wedge\lambda_p.
\endmultline
$$

\proclaim{2.2.2. Proposition} Let $\k$ be a Lie algebra. As a complex of vector
spaces, $\U(\underline{\k})$ is isomorphic to the standard complex $N(\k)$ 
of $\k$.
\qed\endproclaim
In particular, it follows that $N(\k)$ has an algebra structure. The same
structure is described in \cite{MP1, \S 3}. 

We see that our equivariant $(\A,K)$-complexes are DG modules over $N(\k)$.
This suggests the following definitions. Let $\D$ be a DG algebra
with an algebraic action $\chi$ of $K$ by (DG) automorphisms, and with a
morphism $\rho:\k @>>> \D$ of DG Lie algebras, such that the same conditions
as for $\A$ are satisfied: $\rho$ is $K$-equivariant, and the differential of 
the action $\chi$ is the same as the action of $\k$ on $\D$ given by inner 
derivations $[\rho(\xi),-]$. Such a $\D$ is called a Harish-Chandra DG
algebra in \cite{BL2}.
Note that $\rho$ being a morphism of DG Lie algebras implies that its image
is contained in the cycles of $\D^0$.

The main example for $\D$ is the standard complex of $\k$, or of a Lie algebra 
$\g$ such that $(\g,K)$ is a Harish-Chandra pair. Other examples are $\D=\Bbb C$ 
or $\D=\U(\k)$ 
(concentrated in degree 0). In these examples $K$ acts trivially on
$\Bbb C$ and by the adjoint action on $\U(\k)$, $N(\k)$ and $N(\g)$, while 
$\rho$ is the zero map, respectively the natural inclusion.

We now define an $(\A,K,\D)$-module to be a complex $V$ of vector spaces, with 
actions $\pi$ of $\A$, $\nu$ of $K$ and $\omega$ of $\D$,
such that the following conditions hold:
\roster
\item"(AKD1)" $\A$ and $K$ act by chain maps, the $K$-action is algebraic, 
and the $\D$-action is a DG-action;
\item"(AKD2)" $\pi$ and $\omega$ are both $K$-equivariant and they commute;
\item"(AKD3)" The representations $\pi$, $\omega$ and $\nu$ of $\k$ satisfy 
the relation $\pi+\omega=\nu$.
\endroster
A morphism between two $(\A,K,\D)$-modules is a chain map which preserves
all the actions. The category $\M(\A,K,\D)$ of all $(\A,K,\D)$-modules is
clearly abelian.
Examples of interest are the following:

(1) $\D=N(\k)$; then $\M(\A,K,\D)$ is the category of equivariant 
$(\A,K)$-complexes.

(2) $\D=\U(\k)$; then $\M(\A,K,\D)$ is the category of complexes of weak
$(\A,K)$-modules.

(3) $\D=\Bbb C$; then $\M(\A,K,\D)$ is the category of complexes of 
$(\A,K)$-modules.

\noindent We will describe the relationship between these categories in 2.5.2.

Let us now briefly discuss right DG modules.
A right module $M$ over $\D$ is a right DG module if 
$$
d_M(mx)=(d_Mm)x+(-1)^{\deg m}m(d_\D x)
$$
for any homogeneous $x\in\D$ and $m\in M$.

We can also define the opposite DG algebra $\D^{opp}$: as a complex of
vector spaces it is the same as $\D$, and the multiplication is
given by
$$
x\circ y=(-1)^{\deg x \deg y}yx
$$
for homogeneous $x$ and $y$, where $yx$ is the multiplication in $\D$.

Then the category of left (respectively right) DG modules over $\D$ is 
naturally isomorphic to the category of right (respectively left) DG modules 
over $\D^{opp}$. Namely if $M$ is a left DG module
 over $\D$, then it becomes a right DG module over $\D^{opp}$ if we define
$$
m\circ x =(-1)^{\deg x\deg m}xm
$$
for homogeneous $x\in\D$ and $m\in M$. 

In case $\D=\U(\d)$ for a DG Lie algebra $\d$, $\D^{opp}$ can be identified
with $\D$. Namely, let us define $\iota:\d@>>>\U(\d)^{opp}$ by
$$
\iota(X)=^\iota X=-X, \quad X\in\d.
$$
It is easily checked to be a morphism of DG Lie algebras.
Therefore it defines a DG algebra morphism from $\U(\d)$ into $\U(\d)^{opp}$. 
In a Poincar\'e-Birkhoff-Witt basis, $\iota$ acts on monomials by flipping the 
order and putting an appropriate sign. So it is clearly
an isomorphism. We will call $\iota$, understood as a map from $\U(\d)$ into
itself, the principal antiautomorphism of $\U(\d)$. It is a generalization of
the principal antiautomorphism in the case of ordinary Lie algebras. 
It clearly satisfies 
$$
^\iota(uv)=(-1)^{\deg u\deg v}\,\,\,{^\iota v}^\iota u
$$
for homogeneous $u,v\in\U(\d)$. 

Note that we can use $\iota$ to pass from left to right DG modules over 
$\U(\d)$ and vice versa. Indeed, the following lemma holds.

\proclaim{2.2.3. Lemma} If $\omega$ is a left (respectively right) action of 
$\U(\d)$ on $M$, then setting
$$
\tilde{\omega}(u)(m)=(-1)^{\deg u\deg m}\omega(^\iota u)(m)
$$
for homogeneous $u\in\U(\d)$ and $m\in M$ defines a right (respectively left) 
action of $\U(\d)$ on $M$.
\qed\endproclaim

The above principal antiautomorphism is actually the antipode map for a Hopf
algebra structure on $\U(\d)$. We will need this structure in Section 14,
so let us briefly describe it.

Let us first remark that for DG algebras $\D$ and $\E$, we can define their 
tensor product $\D\ootimes\E$. It is equal to $\D\otimes\E$ as a 
vector space, it has the standard tensor product grading and differential (see 
the definition of the tensor algebra above), and the multiplication is given by
$$
(x\otimes y)(x'\otimes y')=(-1)^{\deg y\deg x'}(xx'\otimes yy')
$$
for $x,x'\in\D$ and $y,y'\in\E$. It is straightforward to check that
$\D\ootimes\E$ is a DG algebra.

We now define $c:\d@>>>\U(\d)\ootimes\U(\d)$ by
$$
c(X)=X\otimes 1+1\otimes X,\quad X\in\d.
$$
One easily checks that $c$ is a morphism of DG Lie algebras, hence it
extends in a unique way to a DG morphism $c:\U(\d)@>>>\U(\d)\ootimes\U(\d)$.
This is the coproduct for $\U(\d)$.

The counit $\eps:\U(\d)@>>>\Bbb C$ is defined by the requirements 
$\eps(\d\U(\d))=0$ and $\eps(1)=1$.

It is now straightforward to check the following: 

\proclaim{2.2.4. Proposition} With the above definitions, $\U(\d)$ becomes a 
Hopf algebra. Moreover, the Hopf algebra structure is compatible with the DG 
algebra structure, in the sense that $c$ and $\eps$ are not just algebra 
morphisms, but DG algebra morphisms.
\qed \endproclaim

\subhead 2.3. Homotopic and derived categories of $(\A,K,\D)$-modules
\endsubhead
We first define the translation functor for $(\A,K,\D)$-modules. 
Let $V$ be an $(\A,K,\D)$-module. Then $T(V)=V[1]$ is the same as $V$
as a weak $(\A,K)$-module. The grading is shifted: $T(V)^i=V^{i+1}$, 
the differential is given by $d_{T(V)}=-d_V$, and the $\D$-action is twisted:
$$
\omega_{T(V)}(x)v=(-1)^{\deg x}\omega_V(x)v,
$$
for homogeneous $x\in\D$ and $v\in V$. It is clear that $T(V)$ is again an 
$(\A,K,\D)$-module. If we define the action of $T$ on morphisms just by
shifting the degree: $T(f)^i=f^{i+1}$, then $T$ clearly becomes a functor.
Moreover, $T$ is an autoequivalence of the category $\M(\A,K,\D)$; the inverse 
is given by similar shifting but in the opposite direction.

We say that two $(\A,K,\D)$ morphisms $f$ and $g$ from $V$ to $W$ are homotopic
if there exists $h:V@>>>T^{-1}(W)$, which is an $(\A,K)$-morphism and a 
$\D$-morphism of degree 0, but not necessarily a chain map (so it is not a 
morphism in $\M(\A,K,\D)$), such that 
$$
f-g=hd_V+d_Wh.
$$
For equivariant complexes ($\D=N(\k)$), this condition just means that $f$ and 
$g$ are homotopic as morphisms of complexes of weak $(\A,K)$-modules, via a 
homotopy which anticommutes with $i_\xi$'s.

We now define the homotopic category $K(\A,K,\D)$. Its objects are 
$(\A,K,\D)$-modules, and its morphisms are homotopy classes of morphisms in 
$\M(\A,K,\D)$.

\proclaim{2.3.1. Proposition} The category $K(\A,K,\D)$ is a triangulated 
category.
\endproclaim
\pf The proof is basically the same as for the homotopic category of complexes
over an abelian category. First, it is clear that the above translation functor
descends to $K(\A,K,\D)$; namely, if $h$ is a homotopy from $f$ to 0, then $-h$
shifted in degree by 1 is a homotopy from $T(f)$ to 0.

Next we define the cone of a morphism. Let $f:V@>>>W$ be a morphism in 
$\M(\A,K,\D)$. As a weak $(\A,K)$-module and a graded $\D$-module, the cone
of $f$ is $C_f=T(V)\oplus W$. The differential of $C_f$ is
$$
d_f=\pmatrix  d_{T(V)} &   0 \\
                T(f)   & d_W
    \endpmatrix.
$$
Then one easily checks that $C_f$ is an $(\A,K,\D)$-module. 

The inclusion $i_f:W@>>>C_f$ and the projection $p_f:C_f@>>>T(V)$ are clearly 
$(\A,K,\D)$-morphisms. In this way we obtain the standard triangles
$$
V@>f>>W@>{i_f}>>C_f@>{p_f}>>T(V).
$$
We should show that these triangles are well-defined for a morphism $f$ in
the homotopic category $K(\A,K,\D)$, i.e., that homotopic
$f$ and $g$ give isomorphic triangles (in the obvious sense). Then we can 
define the distinguished triangles as all triangles in $K(\A,K,\D)$ which
are isomorphic (in $K(\A,K,\D)$) to some standard triangle. Finally, we should 
check the axioms of triangulated categories.

Note however that both our homotopies and cones are also homotopies and cones
in the category of complexes of weak $(\A,K)$-modules. They only have
the additional structure of graded $\D$-morphisms, respectively modules.
If we now go through the classical proof of the fact that $K(\M(\A,K)_w)$ is
a triangulated category, we see that for the proof we only need to construct
certain morphisms and homotopies. All these are given by matrices and all the
components of these matrices are clearly graded $\D$-morphisms in our case.
Therefore the classical proof goes through without changes.
\qed\enddemo

To obtain the derived category $D(\A,K,\D)$ of $(\A,K,\D)$-modules, we localize
$K(\A,K,\D)$ with respect to quasiisomorphisms, that is, morphisms in 
$K(\A,K,\D)$ which induce isomorphisms on cohomology. 

\proclaim{2.3.2. Lemma} Quasiisomorphisms in $K(\A,K,\D)$ form a saturated 
localizing class (or multiplicative system) compatible with triangulation.
\endproclaim
\demo
{Proof} The proof of this lemma is analogous to the standard proof of the same 
fact for the homotopic category of an abelian category. The simplest approach 
is to first note that taking $p^{\operatorname{th}}$ cohomology is a
cohomological functor for every $p$, i.e., that for any distinguished triangle
$$
V @>f>> W @>>> Z @>>> T(V) 
$$
in $K(\A,K,\D)$ we have a corresponding long exact sequence of vector spaces
$$
\cdots@>>>H^{p-1}(Z)@>>>H^p(V)@>>>H^p(W)@>>>H^p(Z)@>>>H^{p+1}(V)\cdots.
$$
This immediately implies that $f$ is a quasiisomorphism if and only if
$Z$ is acyclic, i.e., $H^p(Z)=0$ for all $p$. 

The above claim about quasiisomorphisms is therefore equivalent to the fact
that acyclic objects form a null system in $K(\A,K,\D)$ (as defined in 
\cite{KS}, Definition 1.6.6.),
with the additional condition corresponding to the saturation condition:
any direct summand of an object in $\N$ is in $\N$. It is trivial to check
that this is indeed the case. For details, see \cite{P1}, 5.2.2 and 1.1.9.
\qed\enddemo

It is a standard fact (see \cite{KS}, Proposition 1.6.9. or \cite{GM}) 
that if the above conditions are satisfied, then the 
localized category $D(\A,K,\D)$ inherits the structure of a triangulated 
category. The distinguished triangles in $D(\A,K,\D)$ are simply all triangles 
isomorphic in $D(\A,K,\D)$ to distinguished triangles in $K(\A,K,\D)$.
 
\proclaim{2.3.3. Corollary} The category $D(\A,K,\D)$ is a triangulated 
category.
\qed\endproclaim

In case $\D=N(\k)$, which is our main example, the category $D(\A,K,N(\k))$ is
called the {\it equivariant derived category of $(\A,K)$-modules}. We will
sometimes denote it just by $D(\A,K)$, and the corresponding homotopic
category by $K(\A,K)$.

All the above goes through without changes for complexes bounded above, below,
or from both sides. The corresponding derived categories are denoted by
$D^*(\A,K,\D)$, where $*$ can be $-$, $+$, or $b$ (or nothing in the case of the 
full derived category as above).

In case the DG algebra $\D$ is nonpositively graded, i.e., $\D^i=0$ for $i>0$,
as is true in all the examples we wish to consider,
we can also introduce the standard concept of truncations on $D(\A,K,\D)$.

Let $V$ be an $(\A,K,\D)$-module with differential $d$. We define 
$V_{\leq q}$ by
$$
V^p_{\leq q}= \left \{ \matrix V^p & p<q \\
                              \ker d^q & p=q \\
                                0   & p>q.
                      \endmatrix
              \right.
$$
and $V_{\geq q}$ by
$$
V^p_{\geq q}= \left \{ \matrix 0 & p<q \\
                              \coker d^{q-1} & p=q \\
                                V^p   & p>q.
                      \endmatrix
              \right.
$$
Using the fact that $\D$ is nonpositively graded, it is easy to check that 
$V_{\leq q}$ is an $(\A,K,\D)$-submodule of $V$, and that $V_{\geq q}$ is the 
quotient of $V$ by an $(\A,K,\D)$-submodule. In particular, they are both 
$(\A,K,\D)$-modules. Furthermore, the natural projection from $V$ to 
$V_{\geq q+1}$ induces a quasiisomorphism from $V/V_{\leq q}$ to 
$V_{\geq q+1}$. It follows (cf. \cite{KS}, Proposition 1.7.5.) 
that from the short exact sequence
$$
0 @>>> V_{\leq q} @>>> V @>>> V/V_{\leq q} @>>> 0 
$$
in $\M(\A,K,\D)$, we get a distinguished triangle
$$
V_{\leq q} @>>> V @>>> V_{\geq q+1} @>>> T(V_{\leq q}) 
$$
in $D(\A,K,\D)$. Therefore we proved:

\proclaim {2.3.4. Proposition} For any $V$ in $D(\A,K,\D)$ and 
any integer $p$ there is a distinguished triangle
$$
V_{\leq p} @>>> V @>>> V_{\geq p+1} @>>> T(V_{\leq p})
$$
in $D(\A,K,\D)$.
\qed\endproclaim

The standard proof of the following fact now goes through.
Let $D:\M(\A,K)@>>> D^*(\A,K)$ be the functor sending a module $M$ to
the equivariant complex concentrated in degree $0$ and having $M$ as its
zeroth component. Then:

\proclaim{2.3.5. Corollary} The functor $D$ from
$\M(\A,K)$ to $D^b(\A,K)$ is fully faithful. Its image
is a generating class of $D^b(\A,K)$.
\qed\endproclaim

\subhead 2.4. Functors between homotopic categories \endsubhead
In the next sections we will construct various functors on the level of 
$(\A,K,\D)$-modules,
and then we will need to show that they make sense on the level of homotopic
categories. In this section we give a criterion that will be useful in such
situations.

Let us fix two triples, $(\A,K,\D)$ and $(\B,T,\E)$, as in Section 2.2. Let
$F_1$ be an additive functor from $\M(\A,K,\D)$ into $\M(\B,T,\E)$, and let $F_2$ be
an additive functor between the corresponding categories of graded modules, 
$\M^{GR}(\A,K,\D)$ and $\M^{GR}(\B,K,\E)$. Here morphisms of graded modules are 
required to be of degree 0. Assume that the diagram
$$
\CD
\M(\A,K,\D)       @>{F_1}>>   \M(\B,T,\E) \\
 @V{\For}VV                    @V{\For}VV \\
\M^{GR}(\A,K,\D)  @>{F_2}>>   \M^{GR}(\B,T,\E)
\endCD
$$
commutes. Here $\For$ denotes the obvious forgetful functors. Then we will
loosely denote both $F_1$ and $F_2$ by $F$ and think of $F$ as a functor defined 
both on the level of DG modules and on the level of graded modules, in a
compatible way. We also assume that $F$ commutes with translations, i.e., that
both $F_1$ and $F_2$ commute with translations.

Now we need to define the $\Homd$-complex. Let $V$ and $W$ be two objects of 
$\M(\A,K,\D)$. First we define $\Homd_\Bbb C(V,W)$ to be just
the standard $\Homd$ of complexes. In other words,
$$
\Hom^k_\Bbb C(V,W)=\prod_j \Hom_\Bbb C (V^j,W^{j+k})
$$
and the differential is given by $df = d_W\circ f-(-1)^{\deg f}f\circ d_V$
for a homogeneous $f\in \Homd_\Bbb C(V,W)$. More precisely, 
$$
d^kf=(d_W^{j+k}\circ f_j-(-1)^kf_{j+1}\circ d_V^j),
$$
for $f=(f_j)\in\prod_j \Hom_\Bbb C (V^j,W^{j+k})$.

Now we define $\Homd_{(\A,K,\D)}(V,W)$ as the linear subspace of 
$\Homd_\Bbb C(V,W)$, consisting of $(\A,K)$-morphisms whose homogeneous
components $f$ satisfy
$$
f(xv)=(-1)^{\deg f\deg x}xf(v)
$$
for any homogeneous $x\in\D$ and $v\in V$. It is easy to see that 
$\Homd_{(\A,K,\D)}(V,W)$ is actually a subcomplex of $\Homd_\Bbb C(V,W)$,
and that $\Homd_{(\A,K,\D)}$ is an additive bifunctor into the category
$C(Vect_\Bbb C)$ of complexes of vector spaces.

Note that $\Hom_{\M(\A,K,\D)}(V,W)$ consists
of 0-cycles in $\Homd_{(\A,K,\D)}(V,W)$, and that $\Hom_{K(\A,K,\D)}(V,W)$
is the zeroth cohomology space of the complex $\Homd_{(\A,K,\D)}(V,W)$.

Now we return to a functor $F$ as above. 
Since $(\A,K,\D)$-maps of 
degree $k$ between $V$ and $W$ are the same as $(\A,K,\D)$-maps of degree 0 
between $V$ and $W[k]$, and $F$ commutes with translations, we see that $F$ 
defines a linear map
$$
F=F_{V,W}:\Homd_{(\A,K,\D)}(V,W) @>>> \Homd_{(\B,T,\E)}(F(V),F(W)).
$$

\proclaim{2.4.1. Theorem} Let $F$ be a functor as above, such that for any $V$
and $W$ in $\M(\A,K,\D)$ the linear map $F_{V,W}$ defined above is a chain map.
 Then $F$ defines an exact functor from $K(\A,K,\D)$ into $K(\B,T,\E)$.
\endproclaim

\example{Remark} Recall that an exact functor between triangulated categories
is a functor commuting with translation and sending distinguished triangles to 
distinguished triangles.
\endexample

\pf Since $F_{V,W}$ is a chain map, it defines a map on the level of cohomology, 
so it defines a map on the level of homotopic classes of morphisms.
Therefore to finish the proof, it is enough to show that $F$ preserves cones
(we already assumed that $F$ commutes with translation).

Let $f:V@>>>W$ be a morphism in $\M(\A,K,\D)$. Since $F$ is additive and
commutes with translation, $F(C_f)=C_{F(f)}$ as graded $(\B,K,\E)$-modules.
So we only need to show that the differential of $F(C_f)$ is equal to the
differential $d_{F(f)}$ of the cone $C_{F(f)}$, i.e., that
$$
d_{F(C_f)}=\pmatrix  d_{T(F(V))} & 0 \\
                      T(F(f))    & d_{F(W)}
             \endpmatrix.
$$
In the following we denote the inclusions and projections corresponding to
direct sum decompositions by $i$ and $p$ with appropriate indices.

Since $i_W:W@>>>C_f$ and $p_{T(V)}:C_f@>>>T(V)$ are chain maps, they are cycles
in the corresponding $\Homd$ complexes. Thus $F(i_W)=i_{F(W)}$ and 
$F(p_{T(V)})=p_{T(F(V))}$ are also cycles, i.e.,
$$
d_{F(C_f)}i_{F(W)}-i_{F(W)}d_{F(W)}=0
$$
and 
$$
d_{T(F(V))}p_{T(F(V))}-p_{T(F(V))}d_{F(C_f)}=0.
$$
This implies
$$
p_{F(W)}d_{F(C_f)}i_{F(W)}=p_{F(W)}i_{F(W)}d_{F(W)}=d_{F(W)};
$$
$$
p_{T(F(V))}d_{F(C_f)}i_{F(W)}=p_{T(F(V))}i_{F(W)}d_{F(W)}=0;
$$
$$
p_{T(F(V))}d_{F(C_f)}i_{T(F(V))}=d_{T(F(V))}p_{T(F(V))}i_{T(F(V))}=d_{T(F(V))};
$$
so we have checked that three of the matrix coefficients of $d_{F(C_f)}$ agree
with the corresponding matrix coefficients of $d_{F(f)}$. To calculate the
last matrix coefficient, we first calculate the differential in 
$\Homd_{(\A,K,\D)}(C_f,W)$ of the graded morphism $p_W$:
$$
\multline
d(p_W)=d_Wp_W-p_Wd_{C_f}=d_W\pmatrix 0 & 1 \endpmatrix -
     \pmatrix 0 & 1 \endpmatrix \pmatrix d_{T(V)}  &   0  \\
                                          T(f)     &  d_W
                                \endpmatrix
=\\
\pmatrix 0 & d_W \endpmatrix - \pmatrix T(f) & d_W \endpmatrix
=\pmatrix -T(f) & 0 \endpmatrix
=-T(f) p_{T(V)}.
\endmultline
$$
Therefore $F(d(p_W))=d(F(p_W))$ becomes
$$
-F(T(f)p_{T(V)})=d_{F(W)}F(p_W)-F(p_W)d_{F(C_f)},
$$
which implies
$$
p_{F(W)}d_{F(C_f)}=d_{F(W)}p_{F(W)}+T(F(f))p_{T(F(V))}.
$$
Now we can calculate the last matrix coefficient:
$$
p_{F(W)}d_{F(C_f)}i_{T(F(V))}=d_{F(W)}p_{F(W)}i_{T(F(V))}+T(F(f))p_{T(F(V))}i_{T(F(V))}=T(F(f)).
$$
This finishes the proof.
\qed\enddemo

Let us remark that in the classical situation one usually starts with a 
functor on the level of abelian categories, and then lets it act on complexes
by acting  on all components and on differentials. The action on chain
maps is given componentwise. Such a functor obviously
acts on graded objects and morphisms, and it trivially meets the conditions of 
2.4.1. Therefore such a functor always defines a functor on the level of
homotopic categories.

The first functor for which the above conditions are nontrivial is 
$\Homd_{(\A,K,\D)}$ itself. Let $V$ be an $(\A,K,\D)$-module. Then 
$\Homd_{(\A,K,\D)}(V,-)$ is a functor from $\M(\A,K,\D)$ into 
$C(Vect_\Bbb C)$, which is compatibly defined on the level of graded modules; 
in both cases it is defined on morphisms by $\alpha\mapsto \alpha\circ -$. 

It is easy to show that $\Homd_{(\A,K,\D)}(V,-)$ commutes with translation.
It remains to show that the map sending $\alpha\in\Homd_{(\A,K,\D)}(W_1,W_2)$
into 
$$
\tilde\alpha=\alpha\circ - \in \Homd_\Bbb C(\Homd_{(\A,K,\D)}(V,W_1),\Homd_{(\A,K,\D)}(V,W_2))
$$
 is a chain map. Since $\Homd_{(\A,K,\D)}(V,W)$
is a subcomplex of $\Homd_\Bbb C(V,W)$, this follows immediately from 
the first part of the following lemma.

\proclaim{2.4.2. Lemma} Let $V$, $W_1$ and $W_2$ be complexes of vector spaces.
\roster
\item"(i)" For $\alpha\in\Homd_\Bbb C(W_1,W_2)$, define
$\tilde\alpha\in \Homd_\Bbb C(\Homd_\Bbb C (V,W_1),\Homd_\Bbb C (V,W_2))$ by
$\tilde\alpha=\alpha\circ -$. Then the map $\alpha\mapsto\tilde\alpha$ is a chain map.
\item"(ii)" The map $\alpha\mapsto\alpha\otimes id$ from $\Homd_\Bbb C(W_1,W_2)$
into $\Homd_{\Bbb C}(W_1\otimes_\Bbb C V,W_2\otimes_\Bbb C V)$  
is a chain map. \qed
\endroster
\endproclaim

The proof of the lemma is straightforward, and hence we get

\proclaim{2.4.3. Proposition} For any $(\A,K,\D)$-module $V$, 
$\Homd_{(\A,K,\D)}(V,-)$ defines an exact functor from $K(\A,K,\D)$ into 
$K(Vect_\Bbb C)$.
\qed\endproclaim

We remark that similarly one can prove the analogous claim for the first 
variable, and it is then easy to show that $\Homd_{(\A,K,\D)}$ is an exact 
bifunctor from $K(\A,K,\D)^{opp}\times K(\A,K,\D)$ into $K(Vect_\Bbb C)$. 

The last question we want to address in this section is adjunction of
functors on the homotopic level. The following easy fact is enough for
most purposes.

\proclaim{2.4.4. Proposition} Let $F:\M(\A,K,\D)@>>>\M(\B,T,\E)$ be left 
adjoint
to the functor $G:\M(\B,T,\E)@>>>\M(\A,K,\D)$. Suppose that both functors 
define the same
named functors on the level of homotopic categories. Then $F$ and $G$ are 
adjoint functors between the homotopic categories.
\endproclaim

\pf The adjunction morphisms on the level of homotopic categories are defined
as the homotopy classes of the adjunction morphisms on the level of complexes.
The compositions in the definition of adjunction are the identities because
it is so on the level of representatives.
\qed\enddemo

In most situations the reason for $F$ and $G$ to define functors on the level 
of homotopic categories will be the fact that they satisfy the conditions of
2.4.1. If this is the case, and if the adjunction morphisms
$\Phi_X$ and $\Psi_Y$ are natural not only with respect to morphisms (which
are in particular chain maps), but also
with respect to graded morphisms (which are not necessarily chain maps), then 
we actually have the following stronger
version of 2.4.4 (2.4.4 follows from it by taking zeroth cohomology):

\proclaim{2.4.5. Theorem} Let $F:\M(\A,K,\D)@>>>\M(\B,T,\E)$ be left adjoint
to $G:\M(\B,T,\E)@>>>\M(\A,K,\D)$. Assume that both functors satisfy the 
conditions of 2.4.1, and that the adjunction natural transformations
$\Phi$ and $\Psi$ are natural with respect to graded morphisms. Then
$$
\Homd_{(\B,T,\E)}(FX,Y)\cong\Homd_{(\A,K,\D)}(X,GY)
$$
for any $X\in\M(\A,K,\D)$ and $Y\in\M(\B,T,\E)$. 
\endproclaim
\pf 
We will adapt the standard proof of the fact that adjunction of two functors
(defined by the equality $\Hom(FX,Y)=\Hom(X,GY)$)
is equivalent to existence of adjunction natural transformations; see \cite{P2},
A2.1.2.

We define 
$$
\alpha_{X,Y}:\Homd_{(\B,T,\E)}(FX,Y)@>>>\Homd_{(\A,K,\D)}(X,GY)
$$
 and $\beta_{X,Y}$ in the opposite direction by 
$$
\alpha_{X,Y}(f)=G(f)\circ \Phi_X;\quad \beta_{X,Y}(g)=\Psi_Y\circ F(g).
$$
This makes sense since $F$ and $G$ are defined on graded morphisms. The proof
that $\alpha$ and $\beta$ are natural, and that their compositions in both
orders give identity morphisms, is the same as in the classical proof (but
using our present assumptions). The only thing we have to show is that 
$\alpha_{X,Y}$ is a chain map. However, using the fact that $G$ induces a chain
 map between $\Homd$-spaces, and that $\Phi_X$ is a chain map, we see
$$
\multline
\alpha_{X,Y}(df)=G(df)\circ\Phi_X=d(Gf)\circ\Phi_X=
(d_{GY}\circ Gf-(-1)^{\deg f}G(f)\circ d_{GFX})\circ\Phi_X=\\
d_{GY}\circ(Gf\circ\Phi_X)-(-1)^{\deg f}(G(f)\circ\Phi_X)\circ d_X=d(\alpha_{X,Y}(f)).
\endmultline
$$
Namely, $\deg(G(f)\circ\Phi_X)$ is clearly the same as $\deg f$.
\qed\enddemo

\subhead 2.5. Change of DG algebras \endsubhead
Let $(\A,K,\D)$ and $(\A,K,\E)$ be triples as in Section 2.2, with $\phi$, 
$\psi$, $\chi_\E$, $\chi_\D$, $\rho_\E$ and $\rho_\D$ the corresponding maps. 
Let $\alpha:\D@>>>\E$ be a morphism of DG algebras, intertwining the 
$K$-actions, and such that 
$$
\rho_\E=\alpha\circ\rho_\D.
$$
Then $\alpha$ defines a forgetful functor $\For$ from $\M(\A,K,\E)$ into 
$\M(\A,K,\D)$ in the obvious way: the $\D$ action is obtained by composing the 
$\E$-action with $\alpha$. We want to construct adjoints of $\For$.

For simplicity, we assume that both $\D$ and
$\E$ have principal antiautomorphisms, commuting with the $K$-action, which we 
denote by $\iota$ in each case; we also assume that 
$^\iota\rho(\xi)=-\rho(\xi)$ for any $\xi\in\k$ and that the morphism $\alpha$ 
intertwines the $\iota$'s. This will be satisfied in all the examples we
need, and it simplifies matters by allowing us to work with left modules only.

To construct the left adjoint, let $V$ be an $(\A,K,\D)$-module, with
actions $\pi_V$ of $\A$, $\nu_V$ of $K$ and $\omega_V$ of $\D$.
We consider it a right DG module over $\D$ as in 2.2.3. Now we consider
$V\otimes_\Bbb C \E$. It is a right DG module over $\E$ for the right 
multiplication in the second factor; as a complex, it is the standard
tensor product of complexes (see the construction of the tensor algebra 
in 2.2). Let $Z$ be the linear subspace of $V\otimes_\Bbb C \E$ spanned
by all elements of the form
$$
vx\otimes e - v\otimes xe
$$
for (homogeneous) $x\in\D$, $v\in V$ and $e\in\E$. Then $Z$ is easily seen
to be a subcomplex and an $\E$-submodule of $V\otimes_\Bbb C \E$. Therefore,
$V\otimes_\D \E=(V\otimes_\Bbb C \E)/Z$ is a right DG module over $\E$. We
consider it a left $\E$-module by 2.2.3. In other words, the action is given
by
$$
\omega_\E(e)(v\otimes e')=(-1)^{\deg e(\deg v+\deg e')}v\otimes {e'}\,{^\iota e}
$$
for $e,e'\in\E$ and $v\in V$.

We let $\A$ act on $V\otimes_\D\E$ in the first variable, and $K$ in both
variables. In other words, 
$$
\pi(a)=\pi_V(a)\otimes\bold1;\quad \nu(k)=\nu_V(k)\otimes\chi_\E(k)
$$
for $a\in\A$ and $k\in K$. It is straightforward to see that these actions
are well defined, and that in this way $V\otimes_\D\E$ becomes an 
$(\A,K,\E)$-module.

Clearly, in this way we have defined a functor $V\mapsto V\otimes_\D\E$ from
$\M(\A,K,\D)$ into $\M(\A,K,\E)$.
This functor is left adjoint to $\For$. The adjunction morphisms are given
as follows. For an $(\A,K,\D)$-module $V$, $\Phi_V:V@>>>V\otimes_\D\E$ is 
given by 
$$
\Phi_V(v)=v\otimes 1.
$$
For an $(\A,K,\E)$-module $W$, $\Psi_W:W\otimes_\D\E@>>>W$ is given by 
$$
\Psi_W(w\otimes e)=(-1)^{\deg w\deg e}\omega_W(^\iota e)w.
$$
It is easy to see that these morphisms are well defined, that $\Phi_V$ is
an $(\A,K,\D)$-morphism and that $\Psi_W$ is an $(\A,K,\E)$-morphism, and
that in this way we really get adjunction. So we have

\proclaim{2.5.1. Theorem} The functor $V\mapsto V\otimes_\D\E$ described
above is left adjoint to the forgetful functor from $\M(\A,K,\E)$ into 
$\M(\A,K,\D)$. The same is true on the level of homotopic categories.
\endproclaim

\pf It remains to prove the claim about homotopic categories. 
It is clear that $-\otimes_\D \E$ is well-defined on the level of graded 
modules. To check that this functor commutes with translation, we use the 
obvious identification of $T(V)\otimes_\D \E$ and $T(V\otimes_\D \E)$. The 
actions of $\A$ and $K$ obviously agree, and it is straightforward to check 
that the action of $\E$ and the differential agree as well.

To apply 2.4.1, it remains to show that our functor induces chain maps between
corresponding $\Homd$-complexes. However, this follows from 2.4.2.(ii).
\qed\enddemo

\example{Remark} It would look more natural to consider the similarly 
defined functor $V\mapsto \E\otimes_\D V$ instead of the above
functor. This would eliminate passing from left to right modules and
back. On the other hand, in the usual definition of the standard
complex of a Lie algebra the $i_\xi$'s come from the right
multiplication. We do not want to change this, and one of the
important examples will be $\D=\U(\k)$, $\E=N(\k)$.  Hence we consider
$V\otimes_\D\E$, with $\E$-action given by right multiplication.
\endexample

We will not use the other adjoint in this paper, therefore we describe it
only briefly. 
For an $(\A,K,D)$-module $V$, we first consider the complex
$\Homd_\D(\E,V)$. $\Homd_\D$ is defined similarly to $\Homd_{(\A,K,\D)}$
from 2.4; here $\E$ is considered a $\D$-module via left multplication. 
We let $\A$ act on $\Homd_\D(\E,V)$ in the second variable, $K$ on both
variables, and $\E$ by right multiplication in the first variable, i.e., by
$$
(\omega(e)f)(e')=(-1)^{\deg e(\deg f+\deg e')}f(e'e)
$$
for (homogeneous) $f\in\Homd_\D(\E,V)$ and $e,e'\in\E$. 
One easily checks that $\Homd_\D(\E,V)$ satisfies all the properties of an
$(\A,K,\E)$-module, except that the $K$ action is in general not algebraic.
To fix this, we take the $K$-algebraic part, and consider the functor
$$
V\longmapsto \Homd_\D(\E,V)^{alg}.
$$
One now checks that this functor is right adjoint to the forgetful functor.

\example{2.5.2. Examples} Here are the two most important cases of the above
constructions for our applications.
\roster
\item"(i)" $\D=\U(\k)$, $\E=N(\k)$. Then $\M(\A,K,\U(\k))$ is just the
category $C(\M(\A,K)_w)$ of complexes of weak $(\A,K)$-modules, and 
$\M(\A,K,N(\k))$ is the category $C(\A,K)$ of equivariant $(\A,K)$-complexes.
$\For$ is just the obvious forgetful functor, and the left adjoint of $\For$ is
$$
V\longmapsto C_\k(V)=V\otimes_{\U(\k)}N(\k).
$$
Since $N(\k)$ is a free $\U(\k)$-module for the left multiplication, this
functor preserves acyclicity, i.e., if a complex of weak  $(\A,K)$-modules $V$
is acyclic, then $C_\k(V)$ is also acyclic.

$\For$ has also a right adjoint (which we won't need):
$$
C^\k=\Hom_{\U(\k)}(N(\k),-).
$$
Here we do not need to take the algebraic part since $N(\k)$ is finitely 
generated over $\U(\k)$. This functor also preserves acyclicity.

\item"(ii)" $\D=N(\k)$, $\E=\Bbb C$ and $\alpha:N(\k)@>>>\Bbb C$ is the
counit map (see 2.2.4). Then $\M(\A,K,N(\k))=C(\A,K)$ as in (i), while 
$\M(\A,K,\Bbb C)$ is the category
$C(\M(\A,K))$ of complexes of $(\A,K)$-modules. The corresponding forgetful
functor can be described in the following way: a complex $V^\cdot$ of
$(\A,K)$-modules can be viewed as an equivariant $(\A,K)$-complex, if we
define all $i_\xi$'s to be $0$. We will denote this functor by $Q$.
The left adjoint of $Q$ is the functor
$$
V\longmapsto V\otimes_{N(\k)}\Bbb C,
$$
i.e., the functor of taking $N(\k)$-coinvariants. The right adjoint of $Q$
is the functor of taking $N(\k)$-invariants.
\endroster
\endexample

\subhead 2.6. A property of K-projectives \endsubhead
Let us recall that an $(\A,K,\D)$-module $V$ is called K-projective, if
for any acyclic $(\A,K,\D)$-module $W$, $\Hom_{K(\A,K,\D)}(V,W)=0$. Dually,
one defines K-injective $(\A,K,\D)$-modules: $V$ is K-injective if
for any acyclic $(\A,K,\D)$-module $W$, $\Hom_{K(\A,K,\D)}(W,V)=0$.

Since any translate of an acyclic object is acyclic, one sees that $V$ is
K-projective if and only if for any acyclic $(\A,K,\D)$-module $W$, 
the complex $\Homd_{K(\A,K,\D)}(V,W)$ is acyclic. A similar claim holds
for K-injectives. Also, the K-projectives form a null system in $K(\A,K,\D)$,
and so do the K-injectives.
See \cite{P2}, \S\S A1.1, A1.4.

The purpose of this section is to prove a
property of K-projectives, which is actually a property of any null system. 
This property is proved in 2.6.4 below. It says that if an $(\A,K,\D)$-module 
$V$ has a finite filtration such that the corresponding graded modules are 
K-projective, and if this filtration satisfies a certain splitting condition, 
then $V$ is K-projective. We will need this property 
in Section 4.2 to show that the standard complex $N(\k)$ is K-projective 
in a certain category. This will imply that equivariant Zuckerman functors, 
which we define first on the level of homotopic categories, actually make sense
on the level of equivariant derived categories.

The property we need depends on certain properties of the cone, which
seem to be interesting in themselves.
First we define the notion of a semisplit short exact sequence (see \cite{Sp}).
A short exact sequence
$$
0@>>>U@>>>V@>>>W@>>>0
$$
in $\M(\A,K,\D)$ is semisplit if it is split in the category of graded $(\A,K,\D)$-modules.
In other words, we can identify $V$ with $W\oplus U$ as a graded module, but
the inclusion of $W$ is not necessarily a chain map. Clearly, for a morphism
$f:X@>>>Y$ in $\M(\A,K,\D)$, the short exact sequence 
$$
0@>>>Y@>{i_f}>>C_f@>{p_f}>>T(X)@>>>0
$$
is semisplit. We are going to show that any semisplit sequence arises in
this way.

Let $0@>>>U@>>>V@>>>W@>>>0$ be a semisplit short exact sequence. Then the
matrix of the differential of $V$ with respect to the decomposition $V=W\oplus U$ is
$$
d_V=\pmatrix d_W  &  0 \\
            \delta & d_U
    \endpmatrix
$$
where $\delta:W@>>> U$ is a linear map of degree 1. We denote by the same 
letter the corresponding linear map from $W$ into $T(U)$ of degree 0.

\proclaim{2.6.1. Lemma} The map $\delta:W@>>>T(U)$ is a morphism in 
$\M(\A,K,\D)$.
\endproclaim

\pf It is clear that $\delta$ is an $(\A,K)$-morphism. To see it 
is a $\D$-morphism, we first note that for any $x\in\D$,
$$
\omega_V(x)=\pmatrix \omega_W(x) &      0    \\
                          0      & \omega_U(x)
            \endpmatrix
$$
since $V=W\oplus U$ as a $\D$-module. Let 
$v=\pmatrix w \\ u \endpmatrix \in V$.
Since $V$ is a DG module over $\D$, we have
$$
d_V(\omega_V(x)v)=\omega_V(d_\D x)v+(-1)^{\deg x}\omega_V(x)d_V(v),
$$
which becomes 
$$
\pmatrix d_W(\omega_W(x)w) \\
      \delta(\omega_W(x)w)+d_U(\omega_U(x)u)
\endpmatrix
=\pmatrix \omega_W(d_\D x)w+(-1)^{\deg x}\omega_W(x)d_W(w)  \\
       \omega_U(d_\D x)u+(-1)^{\deg x}(\omega_U(x)\delta w +\omega_U(x)d_U u)
\endpmatrix
$$
in the matrix representation. Since $W$ and $U$ are also DG modules over $\D$, 
this implies
$$
\delta\omega_W(x)=(-1)^{\deg x}\omega_U(x)\delta=\omega_{T(U)}(x)\delta,
$$
so $\delta$ is indeed a $\D$-morphism.

It remains to show that $\delta$ is a chain map. However, $d_V^2=0$, i.e.,
$$
\pmatrix        d_W^2        &    0 \\
        \delta d_W+d_U\delta &  d_U^2
\endpmatrix
=0,
$$
hence
$$
\delta d_W=-d_U \delta =d_{T(U)}\delta.
$$
This finishes the proof.
\qed\enddemo

Now we see that $V$ is the cone of the morphism 
$T^{-1}(\delta):T^{-1}(W)@>>>U$. Namely, $V=T(T^{-1}(W))\oplus U$ as a 
graded module, and its differential
is $\pmatrix d_{T(T^{-1}W)}  &   0  \\
             T(T^{-1}\delta) &  d_U
    \endpmatrix$,
which is equal to $d_{T^{-1}\delta}$. Therefore we have proved the following
characterization of the cone.

\proclaim{2.6.2. Proposition} Let $0@>>>U@>>>V@>>>W@>>>0$ be a semisplit short 
exact sequence in $\M(\A,K,\D)$. Then any choice of a splitting for this 
sequence exhibits $V$ as the cone of a morphism from $T^{-1}(W)$ into $U$.
\qed\endproclaim

This has the following immediate consequences, which are useful for proving 
that certain objects are K-projective or K-injective (see 3.2.6).

\proclaim{2.6.3. Corollary} Let $0@>>>U@>>>V@>>>W@>>>0$ be a semisplit short 
exact sequence in $\M(\A,K,\D)$. Let $\N$ be a null system in $K(\A,K,\D)$.
Then if any two of the objects $U$, $V$ and $W$
are in $\N$, so is the third.
\endproclaim
\pf From 2.6.2 we see that we have a distinguished triangle
$$
T^{-1}(W)@>>>U@>>>V@>>>W
$$
in $K(\A,K,\D)$. Therefore the claim follows from the properties (N2) and (N3)
of a null system.
\qed\enddemo

\proclaim{2.6.4. Corollary} Let $\N$ be a null system in $K(\A,K,\D)$ and let
$V$ be an $(\A,K,\D)$-module. Assume
$$
0=F_0V\subset F_1V\subset\dots\subset F_nV=V
$$
is a finite filtration of $V$ by $(\A,K,\D)$-submodules, such that the
graded objects
$$
Gr_iV=F_iV/F_{i-1}V,\quad i=1,\dots,n
$$
belong to $\N$. Assume further that $F_iV$ is isomorphic to 
$F_{i-1}V\oplus Gr_iV$ as a graded $(\A,K,\D)$-module. 
Then $V$ belongs to $\N$.
\endproclaim
\pf Starting from $0\in\N$ and using 2.6.3, we prove inductively that $F_iV$ 
belong to $\N$ for all $i$. In particular, for $i=n$ we get that $V$ belongs 
to $\N$.
\qed\enddemo

It is possible to generalize 2.6.4 to infinite filtrations, both increasing
and decreasing (for decreasing filtrations, one needs a finiteness assumption).
These generalizations together with the results of 2.5 lead to constructions of
enough K-injectives in $K(\A,K,\D)$, and also K-projectives in case $K$ is
reductive. The construction of K-projectives is explained in \cite{P2}, \S 1,
while the construction of K-injectives was done in \cite{BL2}, 1.15.3, in 
essentially the same way. Both constructions are explained in detail in 
\cite{P1}, \S 5.6.

\head  3. Equivariant Zuckerman  functors
\endhead

In this section we consider a situation as in Section 1.3: $(\A,K)$ is a 
Harish-Chandra pair, and $\gamma:T@>>>K$ is a morphism of algebraic groups.
We will assume that the differential of $\gamma$ is injective, so that the Lie 
algebra $\t$ of $T$ can be viewed as a Lie subalgebra of the Lie algebra $\k$ 
of $K$, and $N(\t)$ can be viewed as a DG subalgebra of $N(\k)$. Note that 
$(\U(\k),T)$ is another Harish-Chandra pair in the obvious way.

In this section we are going to construct the right adjoint $\Gamma^{eq}_{K,T}$
of the forgetful functor from the category of equivariant $(\A,K)$-complexes 
(i.e., $(\A,K,N(\k))$-modules) into the category of equivariant 
$(\A,K)$-complexes (i.e., $(\A,T,N(\t))$-modules). $\Gamma^{eq}_{K,T}$ will  
also define a functor between the
corresponding homotopic categories. Furthermore, we will show that, in case 
$T$ is reductive, this functor is ``acyclic'' (i.e., preserves acyclic 
complexes, or equivalently, preserves quasiisomorphisms). It follows that this 
functor defines a functor between the corresponding equivariant derived 
categories, denoted again by $\Gamma^{eq}_{K,T}$.
This is the equivariant analogue of the derived Zuckerman functor. 
Finally, $\Gamma^{eq}_{K,T}$ is given on all levels by the same explicit
formula (not involving resolutions), analogous to the Duflo-Vergne formula
from 1.3.3.

We will further show that in case $\A$ is a flat $\U(\k)$-module for the 
right multiplication, the cohomology modules of $\Gamma^{eq}_{K,T}(V)$, for an
$(\A,T)$-module $V$ viewed as an equivariant complex concentrated in degree 0,
are equal to the classical derived Zuckerman functors of $V$.  

One advantage of the equivariant construction over the classical one
is the fact that the equivariant construction is independent of $\A$; in 
particular, if
$(\g,K)$ is a Harish-Chandra pair, then the same formula describes the adjoint 
for $\A=\U(\g)$ and for $\A=\U_\theta$, the quotient of $\U(\g)$ corresponding
to an infinitesimal character. This leads to a natural way of localizing the
Zuckerman construction; see \cite{MP1}, \S 4, and \cite{MP2}.

\subhead 3.1. Construction of equivariant Zuckerman functors \endsubhead
Let $V$ be an equivariant $(\A,T)$-complex, i.e., an $(\A,T,N(\t))$-module. 
Let $\pi_V$, $\nu_V$ and $\omega_V$ be the corresponding actions. 
We consider $V$ as a complex of weak $(\A,T)$-modules and apply to it the
functor $\Ind_w$ from Section 1.3. In this way we get a complex of weak 
$(\A,K)$-modules $\Ind_w(V)=R(K,V)=R(K)\otimes_\Bbb C V$. We denote the 
corresponding actions of $\A$ and $K$ by $\pi_R$ and $\nu_R$; recall that 
$\nu_R$ is just the right translation, and that for $a\in\A$
$$
(\pi_R(a)F)(k)=\pi_V(\phi_K(k)a)(F(k)),\quad k\in K.
$$
Both of these are actions by chain maps. Note that the differential of $R(K,V)$
is induced by the differential of $V$.

Now we want to construct a $(\k,T,N(\t))$-action on $R(K,V)$. We
already constructed actions $\lambda_\k$ of $\k$ and $\lambda_T$
of $T$ in Section 1.3; they were both given as the tensor product of
the left regular action on $R(K)$ with the given action on $V$.
Clearly, in our present situation when $V$ is a complex, they are
actions by chain maps (namely, the actions $\pi_V$ and $\nu_V$ are by
chain maps). Therefore, $R(K,V)$ is a complex of weak
$(\k,T)$-modules with respect to $\lambda_\k$ and $\lambda_T$.

We define an action of $N(\t)$ on $R(K,V)$ by
$$
(\lambda_N(n)F)(k)=\omega_V(n)(F(k)),\quad k\in K,
$$
for $n\in N(\t)$ and $F\in R(K,V)$. In other words, $\lambda_N$ is the tensor 
product of the trivial action of $N(\t)$ on $R(K)$ with the given action on 
$V$. It is clearly a DG action.

The following lemma is straightforward.

\proclaim{3.1.1. Lemma} Let $V$ be an $(\A,T,N(\t))$-module. Then
\roster
\item"(i)" With the above described action $\lambda$, $\Ind_w(V)=R(K,V)$ is a
$(\k,T,N(\t))$-module; 
\item"(ii)" The actions $\pi_R$ and $\nu_R$ of $\A$ and $K$ on $R(K,V)$ commute
with the $\lambda$-action;
\item"(iii)" Let $\varphi:V@>>>W$ be a morphism of $(\A,T,N(\t))$-modules.
Since $\varphi$ is a morphism of complexes of weak $(\A,T)$-modules, we
can consider the morphism $\Ind_w(\varphi)$ of complexes of weak 
$(\A,K)$-modules. Then $\Ind_w(\varphi)$ intertwines the $\lambda$-actions. \qed
\endroster
\endproclaim

Consider now the standard complex $N(\k)$ of $\k$ as a
$(\k,T,N(\t))$-module, in the following way.  $\k$ acts by
$\pi_N$ which is the left multiplication, $T$ acts by the action
$\nu_N$ induced by $\phi$, and $N(\t)$ acts by $\omega_N$ which is the
right multiplication twisted to a left action: 
$$
\omega_N(m)n=(-1)^{\deg m\deg n}n^\iota m
$$
for $m\in N(\t)$ and $n\in N(\k)$. Here $\iota$ is the principal 
antiautomorphism of $N(\t)$, which is the same as the restriction 
of the principal antiautomorphism of $N(\k)$. See Section 2.2, in particular
2.2.3.

It is straightforward to check that in this way $N(\k)$ indeed becomes 
a $(\k,T,N(\t))$-module. We define
$$
\Gamma^{eq}_{K,T}(V)=\Homd_{(\k,T,N(\t))}(N(\k),R(K,V)),
$$
where the $(\k,T,N(\t))$-action on 
$R(K,V)$ is the $\lambda$-action. $\Gamma^{eq}_{K,T}(V)$ is a complex of 
vector spaces as we saw at the beginning of Section 2.4.
 
Recall from Section 2.4 that the condition for $f:N(\k)@>>>R(K,V)$ to be an 
$N(\t)$-morphism is
$$
f(\omega_N(m)n)=(-1)^{\deg f\deg m}\lambda_N(m)(f(n)).
$$

We define actions of $\A$, $K$ and $N(\k)$ on $f\in\Gamma^{eq}_{K,T}(V)$ as follows:
$$
\pi_\Gamma(a)f=\pi_R(a)\circ f,\quad a\in\A,
$$
$$
\nu_\Gamma(k)f=\nu_R(k)\circ f,\quad k\in K,
$$
$$
(\omega_\Gamma(n)f)(m)(k)=(-1)^{\deg n\deg f}f(\Ad k(^\iota n)m)(k)
$$
for $a\in\A$, $m,n\in N(\k)$ and $k\in K$.

\proclaim{3.1.2. Theorem} With the above actions, $\Gamma^{eq}_{K,T}(V)$ is an 
$(\A,K,N(\k))$-module. Furthermore, $\Gamma^{eq}_{K,T}$ is a functor from 
$\M(\A,T,N(\t))$ into $\M(\A,K,N(\k))$.
\endproclaim
\pf It is clear from 3.1.1 that $\pi_\Gamma$ and $\nu_\Gamma$ are well-defined 
actions by chain maps. It is also clear that $\nu_\Gamma$ is algebraic, and that 
$\pi_\Gamma$ is $\nu_\Gamma$-equivariant.

Let us prove that the action $\omega_\Gamma$ is a well-defined DG action. 
Let $f\in\Gamma^{eq}_{K,T}(V)$ and $n\in N(\k)$. First, it is obvious that 
$\omega_\Gamma(n)f$ is a linear map of
degree $\deg n+\deg f$. We have to prove that $\omega_\Gamma(n)f$ is a 
$(\k,T,N(\t))$-morphism; then it will be clear that $\omega_\Gamma$ is a DG 
action. Namely, it is clear that $(fn)(m)=f(nm)$, $n,m\in N(\k)$, 
$f\in \Gamma^{eq}_{K,T}(V)$, defines a right DG action of 
$N(\k)$. Hence the same is true for
$$
(fn)(m)(k)=f(\Ad k(n)m)(k),\quad k\in K,
$$
and so, by 2.2.3,  $\omega_\Gamma$ is a left DG action of $N(\k)$.

It is straightforward to check that $\omega_\Gamma(n)f$ is a $T$-morphism and
an $N(\t)$-morphism. 
To see that it is also a $\k$-morphism, let $X\in\k$. Then
$$
\multline
\lambda(X)((\omega_\Gamma(n)f)(m))(k)=\pi_V(X)((\omega_\Gamma(n)f)(m)(k))+
L_X((\omega_\Gamma(n)f)(m))(k)=\\
(-1)^{\deg n\deg f}\pi_V(X)(f(\Ad k(^\iota n)m)(k))+
{d\over dt}((\omega_\Gamma(n)f)(m))(\exp(-tX)k)\big|_{t=0}=\\
(-1)^{\deg n\deg f}(\pi_V(X)(f(\Ad k(^\iota n)m)(k))+
{d\over dt}f(\Ad(\exp(-tX)k)(^\iota n)m)(\exp(-tX)k)\big|_{t=0})=\\
(-1)^{\deg n\deg f}(\pi_V(X)(f(\Ad k(^\iota n)m)(k))+\\
f(-[X,\Ad k(^\iota n)]m)(k)+
{d\over dt}f(\Ad k(^\iota n)m)(\exp(-tX)k)\big|_{t=0})=\\
(-1)^{\deg n\deg f}(f(X\Ad k(^\iota n)m)(k)-f([X,\Ad k(^\iota n)]m)(k))=\\
(-1)^{\deg n\deg f}f(\Ad k(^\iota n)Xm)(k)=(\omega_\Gamma(n)f)(Xm)(k),
\endmultline
$$
for any $m\in N(\k)$ and $k\in K$. Here we used the fact that since $f$ is a 
$\k$-morphism, 
$$
\pi_V(X)(f(\Ad k(^\iota n)m)(k))+
{d\over dt}f(\Ad k(^\iota n)m)(\exp(-tX)k)\big|_{t=0}=f(X\Ad k(^\iota n)m)(k),
$$
for any (fixed) $k\in K$.

The $\nu_\Gamma$-equivariance of $\omega_\Gamma$ and the commuting of
$\omega_\Gamma$ and $\pi_\Gamma$ are straightforward.

It remains to show that $\nu_\Gamma-\pi_\Gamma=\omega_\Gamma$ on $\k$. 
This follows from 1.3.1, which says that for any $F\in R(K,V)$ and $X\in\k$
we have
$$
(\nu_R(X)F)(k)-(\pi_R(X)F)(k)=(\omega_R(X)F)(k)=-(\lambda_\k(\Ad k(X))F)(k),
$$
for any $k\in K$. 

This shows that $\Gamma^{eq}_{K,T}(V)$ is indeed an $(\A,K,N(\k))$-module. 
It is now easy to check that $\Gamma^{eq}_{K,T}$ is a functor; its action on a
morphism $\varphi:V@>>>W$ of $(\A,T,N(\t))$-modules is given by
$$
(\Gamma^{eq}_{K,T}(\varphi)(f))(n)(k)=\varphi(f(n)(k)),
$$
for $f\in\Gamma^{eq}_{K,T}(V)$, $n\in N(\k)$ and $k\in K$.
\qed\enddemo

\proclaim{3.1.3. Theorem} The functor $\Gamma^{eq}_{K,T}$ is right adjoint to
the forgetful functor from the category $\M(\A,K,N(\k))$ into the category
$\M(\A,T,N(\t))$.
\endproclaim
\pf
We will just define the adjunction morphisms and leave the straightforward
checking to the reader. For an $(\A,K,N(\k))$-module $V$,
we define $\Phi_V:V@>>>\Gamma^{eq}_{K,T}(V)$ by
$$
(\Phi_V(v)(n))(k)=(-1)^{\deg v\deg n}\omega_V(^\iota n)\nu_V(k)v,
$$
for $v\in V$, $n\in N(\k)$ and $k\in K$. For an $(\A,T,N(\t))$-module $W$ we 
define $\Psi_W:\Gamma^{eq}_{K,T}(W)@>>>W$ by
$$
\Psi_W(f)=f(1)(1)
$$
for $f\in\Gamma^{eq}_{K,T}(W)$. 
\qed\enddemo

Now we want to show that $\Gamma^{eq}_{K,T}$ defines an exact functor on the 
level of homotopic categories. The proof of this is essentially the same as 
the proof of 2.4.3. Namely, the proof of functoriality of 
$\Gamma^{eq}_{K,T}$ in 3.1.2 shows also that $\Gamma^{eq}_{K,T}$ transforms 
graded $(\A,T,N(\t))$-morphisms into graded $(\A,K,N(\k))$-morphisms. On the 
other hand, it was proved in the proof of 2.4.3 that $\Gamma^{eq}_{K,T}$ 
satisfies the other conditions of 2.4.1; namely, these other conditions are 
independent of the $(\A,K,N(\k))$-action and depend only on the structure of 
a complex of vector spaces. 

So $\Gamma^{eq}_{K,T}$ indeed defines an exact functor between the
homotopic categories, which we denote again by $\Gamma^{eq}_{K,T}$. By 
2.4.4 (or 2.4.5), it is still adjoint to the forgetful functor on the 
homotopic level. Hence we have proved

\proclaim{3.1.4. Proposition} $\Gamma^{eq}_{K,T}$ defines an exact functor 
from the category $K(\A,T,N(\t))$ into the category $K(\A,K,N(\k))$. This 
functor is right adjoint to the forgetful functor.
\qed\endproclaim

To pass to derived categories, we assume that $T$ is reductive (for the case
of non-reductive $T$, see \cite{MP1}, \S 2). In that case, we are going to prove 
that $\Gamma^{eq}_{K,T}$ is an acyclic functor, i.e., transforms
acyclic equivariant $(\A,T)$-complexes into acyclic equivariant 
$(\A,K)$-complexes. This implies (see \cite{P2}, \S A1.2) that it also transforms 
quasiisomorphisms into quasiisomorphisms, and hence defines a functor on the 
level of derived categories; this functor is equal to 
$\Gamma^{eq}_{K,T}$ on objects, and acts on triples in the obvious way.

To show this, note first that if $V$ is an acyclic equivariant 
$(\A,T)$-complex, then $R(K,V)$ is also acyclic, either as a complex of weak 
$(\A,K)$-modules, or as a $(\k,T,N(\t))$-module. On the other hand, 
acyclicity of $\Gamma^{eq}_{K,T}(V)$ can be checked on the level of
vector spaces. Therefore, it is
enough to show that the functor $\Homd_{(\k,T,N(\t))}(N(\k),-)$ transforms
acyclic $(\k,T,N(\t))$-modules into acyclic complexes of vector spaces.

In other words, we need to show that $N(\k)$ is a K-projective 
$(\k,T,N(\t))$-module (see \cite{P2}, \S A1.4 and also \S 2.6 in this paper). 
This will be proved in 3.2.6 
below. Since $\Gamma^{eq}_{K,T}(V):D(\A,T,N(\t))@>>>D(\A,K,N(\k))$ remains 
adjoint to the forgetful functor (see \cite{P2}, A2.3.2), this will imply the 
following theorem.

\proclaim{3.1.5. Theorem} Assume $T$ is reductive. Then the functor 
$\Gamma^{eq}_{K,T}$ is acyclic and
hence defines the same named functor on the level of derived categories,
from $D(\A,T,N(\t))$ into $D(\A,K,N(\k))$.
This functor is right adjoint to the forgetful functor.
\qed\endproclaim

\subhead 3.2. K-projectivity of $N(\k)$ \endsubhead
To prove that $N(\k)$ is a K-projective $(\k,T,N(\t))$-module, we use
the Hochschild-Serre filtration (see \cite{HS}), which we now describe.
To construct this filtration we do not need $T$ to be reductive; this will be 
required only in the proof of K-projectivity.

Let $p$ be an integer between 0 and $\dim\k$. We consider 
$\oplus_{i=0}^p N^{-i}(\k)$. This is clearly a subcomplex and a 
$(\k,T)$-submodule of $N(\k)$; however
it is not an $N(\t)$-submodule. Therefore we define
$$
F_pN(\k)=\omega_N(N(\t))(\oplus_{i=0}^p N^{-i}(\k)).
$$
This is now clearly an $N(\t)$-submodule; however it is also a $\k$-submodule 
since $\pi_N$ and $\omega_N$ commute, and a $T$-submodule since $\omega_N$ is 
$T$-equivariant and $N(\t)$ is $T$-invariant. Finally it is a subcomplex by the 
DG property of $\omega_N$. Hence we have defined an increasing filtration of 
$N(\k)$ by $(\k,T,N(\t))$-submodules. 

In the following lemma we describe $F_pN(\k)$ more explicitly. In particular,
from this description it will become clear that our filtration is really the
Hochschild-Serre filtration. The proof is straightforward.

\proclaim{3.2.1. Lemma} For $k\leq p$, $F_pN^{-k}(\k)=N^{-k}(\k)$. For 
$k\geq p$,
$$
F_pN^{-k}(\k)=N^{-p}(\k){\tsize\bigwedge}^{k-p}(\t).
$$
Here the right hand side is a subspace of $N^{-k}(\k)$ via multiplication in 
$N(\k)$.
\qed\endproclaim

In this description of the filtration, the containment 
$F_pN(\k)\subset F_{p+1}N(\k)$ becomes
$$
N^{-p}(\k){\tsize\bigwedge}^{k-p}(\t)\subset N^{-p-1}(\k)
{\tsize\bigwedge}^{k-p-1}(\t),
$$
the inclusion being given in the obvious way, i.e., by
$$
(u\otimes\xi)(\eta_1\wedge\dots\wedge\eta_{k-p})=
(u\otimes\xi\wedge\eta_1)(\eta_2\wedge\dots\wedge\eta_{k-p}).
$$
Furthermore, it is clear that $F_{\dim\k}N(\k)=N(\k)$, while 
$F_0N(\k)=\U(\k) {\tsize\bigwedge}(\t)$ can be 
identified with $\U(\k)\otimes_\Bbb C{\tsize\bigwedge}(\t)$. We define 
$F_{-1}N(\k)=0$.

We want to identify the graded object corresponding to this filtration.
For this we need the following lemma about exterior algebras, the proof of
which follows easily from the results in \cite{Bo}, III.7.

\proclaim{3.2.2. Lemma} Let $\k$ be a finite dimensional vector space and
$\t$ a subspace of $\k$. Denote by $\pi$ the canonical graded map from 
${\tsize\bigwedge}(\k)$
into ${\tsize\bigwedge}(\k/\t)$ induced by the projection $\k@>>>\k/\t$. Then
\roster
\item"(i)" $\pi$ is surjective and its kernel is ${\tsize\bigwedge}(\k)\t$, 
the ideal in ${\tsize\bigwedge}(\k)$ generated by $\t$;
\item"(ii)" The bigraded linear map 
$$
m:{\tsize\bigwedge}(\k)\otimes_\Bbb C{\tsize\bigwedge}(\t)@>>>{\tsize\bigwedge}
(\k){\tsize\bigwedge}(\t)
$$
given by multiplication is surjective and its kernel is contained in
${\tsize\bigwedge}(\k)\t\otimes_\Bbb C{\tsize\bigwedge}(\t)$;
\item"(iii)" The map 
$\pi^i\otimes 1:{\tsize\bigwedge}^i(\k)\otimes_\Bbb C{\tsize\bigwedge}^j(\t)
@>>>{\tsize\bigwedge}^i(\k/\t)\otimes_\Bbb C{\tsize\bigwedge}^j(\t)$ defines 
a map 
$$
\varphi:{\tsize\bigwedge}^i(\k){\tsize\bigwedge}^j(\t)@>>>{\tsize\bigwedge}^i
(\k/\t)\otimes_\Bbb C{\tsize\bigwedge}^j(\t)
$$
such that $\pi^i\otimes 1=\varphi\circ m$. The kernel of $\varphi$ is
$$
{\tsize\bigwedge}^{i-1}(\k){\tsize\bigwedge}^{j+1}(\t)=({\tsize\bigwedge}^{i-1}
(\k)\t){\tsize\bigwedge}^j(\t)\subset{\tsize\bigwedge}^i(\k){\tsize\bigwedge}^j
(\t). \qed
$$
\endroster
\endproclaim

\example{3.2.3. Remark} We are, of course, going to apply 3.2.2 in the case 
studied above, of a Lie algebra $\k$ and its subalgebra $\t$. In that case, 
$T$ acts on $\k$, $\t$ and $\k/\t$, and hence on the corresponding exterior 
algebras. It is then obvious that all the maps considered in 3.2.2 are 
$T$-morphisms. It is also clear that the maps $m$ and $\varphi$ from 3.2.2 
are morphisms with respect to the action of ${\tsize\bigwedge}(\t)$
via right multiplication.
\endexample

We are ready now to describe the graded object corresponding to our 
filtration. Consider the algebraic $T$-module ${\tsize\bigwedge}^p(\k/\t)$ 
as a weak $(\Bbb C,T)$-module. We can apply the functor $\ind_{\U(\k),\Bbb C}$
from Section 1.2 to this module. In this way we get a weak 
$(\k,T)$-module
$$
\ind_{\U(\k),\Bbb C}{\tsize\bigwedge}^p(\k/\t)=
\U(\k)\otimes_\Bbb C {\tsize\bigwedge}^p(\k/\t);
$$
recall that $\k$ acts on this module via left multiplication in the first 
factor, while $T$ acts on both factors. Now we can consider 
$\U(\k)\otimes_\Bbb C {\tsize\bigwedge}^p(\k/\t)[p]$;
it is a complex of weak $(\k,T)$-modules concentrated in degree $-p$. We now 
apply the functor $C_\t$ from 2.5.2.(i) to this complex. In this way we obtain 
a $(\k,T,N(\t))$-module
$$
C_\t(\ind_{\U(\k),\Bbb C}{\tsize\bigwedge}^p(\k/\t)[p])=
(\U(\k)\otimes_\Bbb C{\tsize\bigwedge}^p(\k/\t)[p])\otimes_{\U(\t)}N(\t).
$$ 
Recall from \S 2.5 that here the tensor product is taken with respect to the 
right action 
of $\U(\t)$ on the first factor, which is the $\omega$-action twisted by 
$\iota$, and the left multiplication in the second factor. The action of 
$\k$ on the above module is given by the left multiplication in the first 
factor, $T$ acts on all three factors, while $N(\t)$ acts on the third factor 
by right multiplication twisted by $\iota$. The differential is the tensor 
product differential, but since 
$\U(\k)\otimes_\Bbb C{\tsize\bigwedge}^p(\k/\t)[p]$ has zero differential, 
we see that
$$
d=(-1)^p\,\,1\otimes d_N,
$$
where $d_N$ denotes the differential of $N(\t)$.

Finally, let us note that that the above complex is concentrated in degrees
between $-p$ and $-p-\dim\t$.

\proclaim{3.2.4. Proposition} The graded object attached to the 
Hochschild-Serre filtration $F\,N(\k)$ is given by the above described modules:
$$
Gr_p N(\k)=C_\t(\ind_{\U(\k),\Bbb C}{\tsize\bigwedge}^p(\k/\t)[p]).
$$
\endproclaim
\pf
Let us first note that for $p\leq k\leq p+\dim\t$, 
$C_\t(\ind_{\U(\k),\Bbb C}{\tsize\bigwedge}^p(\k/\t)[p])^{-k}$ can
be identified with 
$$
\U(\k)\otimes_\Bbb C {\tsize\bigwedge}^p(\k/\t)\otimes_\Bbb C 
{\tsize\bigwedge}^{k-p}(\t).
$$
In this interpretation, the actions of $\k$, $T$ and $\underline{\t}^{-1}$
are given in the obvious way (as before), while the action of 
$X\in\underline{\t}^0$
is still by right multiplication, but we have to use the commuting rules as
follows:
$$
\omega(X)(u\otimes \lambda\otimes\mu)=
-uX\otimes\lambda\otimes\mu+u\otimes[X,\lambda]\otimes\mu+
u\otimes\lambda\otimes[X,\mu].
$$
Let us define a graded linear map 
$f:F_pN(\k)@>>>C_\t(\ind_{\U(\k),\Bbb C}{\tsize\bigwedge}^p(\k/\t)[p])$ as 
follows: for $k<p$, $f^{-k}=0$, while for $k\geq p$,
$$
f^{-k}:\U(\k)\otimes_\Bbb C {\tsize\bigwedge}^p(\k){\tsize\bigwedge}^{k-p}(\t)
@>>>\U(\k)\otimes_\Bbb C {\tsize\bigwedge}^p(\k/\t)\otimes_\Bbb C 
{\tsize\bigwedge}^{k-p}(\t)
$$
is $1\otimes\varphi$, where 
$\varphi:{\tsize\bigwedge}^p(\k){\tsize\bigwedge}^{k-p}(\t)@>>>
{\tsize\bigwedge}^p(\k/\t)\otimes_\Bbb C {\tsize\bigwedge}^{k-p}(\t)$ is 
the map from 3.2.2.(iii). Here we used
the description of $F_pN(\k)$ given by 3.2.1.

By 3.2.2.(iii), $f$ is surjective and its kernel is equal to
$$
\U(\k)\otimes_\Bbb C {\tsize\bigwedge}^{p-1}(\k){\tsize\bigwedge}^{k-p+1}(\t)
$$
in degrees $-k$, for $k\geq p$. It is clear that for $k<p$, 
$\ker f^{-k}=N^{-k}(\k)$. So, by 3.2.1, the kernel of $f$ is exactly 
$F_{p-1}N(\k)$. Therefore to prove the theorem, it only remains to show 
that $f$ is a morphism of $(\k,T,N(\t))$-modules.

It is obvious that $f$ intertwines the $\k$-actions. By 3.2.3 it also 
intertwines the $T$-actions and the actions of $\underline{\t}^{-1}$. To see 
that it also intertwines the actions of $\underline{\t}^0$, note that these 
actions are given on both modules by right multiplication on $\U(\k)$ tensored 
with the adjoint action on the other factor. However, the map $\varphi$ 
intertwines the adjoint actions of $T$, hence also of $\t$. 

In particular, we see that if $k\geq p$, then for $u\otimes\lambda\in N^{-p}(\k)$
and $n\in N^{-k+p}(\t)$, we have
$$
f^{-k}((u\otimes\lambda)\cdot n)=f^{-p}(u\otimes\lambda)\cdot n,
$$
where $-\cdot n$ denotes the right multiplication by $n$ on $N(\k)$, respectively
on the last factor $N(\t)$ in 
$(\U(\k)\otimes_\Bbb C{\tsize\bigwedge}^p(\k/\t)[p])\otimes_{\U(\t)}N(\t)$.
Note also that $(u\otimes\lambda)\cdot n\in F_pN^{-k}(\k)$, so $f^{-k}$ can be
applied to it. Using this, let us prove that $f$ is a chain map.

Let us fix $k>p$; for $k\leq p$ there is nothing to prove. Let
$(u\otimes\xi)\eta$ be an element of 
$F_pN^{-k-1}(\k)=N^{-p}(\k){\tsize\bigwedge}^{k+1-p}(\t)$. Then
$$
f^{-k}(d_N((u\otimes\xi)\eta))=f^{-k}(d_N(u\otimes\xi)\eta +
(-1)^p (u\otimes\xi)d_N(\eta))
$$
by the properties of $d_N$. However, $d_N(u\otimes\xi)\eta$ is in 
$F_{p-1}N^{-k}(\k)$ since $\eta\in{\tsize\bigwedge}^{k-p+1}(\t)$, 
so $f^{-k}$ sends it to 0. Hence the above expression is equal to
$$
(-1)^p f^{-k}((u\otimes\xi)d_N(\eta))=(-1)^p f^{-p}(u\otimes\xi)\cdot d_N(\eta);
$$
here we used the above remark. Now by 3.2.2.(iii),
$$
f^{-p}(u\otimes\xi)=u\otimes\varphi(\xi)=u\otimes\varphi(\xi\cdot 1)=
u\otimes\pi^p(\xi)\otimes 1,
$$
where 
$\pi^p:{\tsize\bigwedge}^p(\k)@>>>{\tsize\bigwedge}^p(\k/\t)$ is induced by the 
projection $\k@>>>\k/\t$. So the above expression is equal to
$$
(-1)^p u\otimes\pi^p(\xi)\otimes d_N(\eta)=d(u\otimes \pi^p(\xi)\otimes\eta)=
d(u\otimes\varphi(\xi\eta))=d(f^{-k-1}((u\otimes\xi)\eta)).
$$
So $f$ is a chain map and this finishes the proof.
\qed\enddemo

Assume now that $T$ is reductive. Then ${\tsize\bigwedge}^p(\k/\t)$ is
a projective $T$-module, and hence $\ind_{\U(\k),\Bbb
C}{\tsize\bigwedge}^p(\k/\t)$ is a projective weak $(\k,T)$-module
since $\ind_{\U(\k),\Bbb C}$ is left adjoint to the forgetful functor
from $\M(\k,T)_w$ into $\M(T)$.  This implies that
$\ind_{\U(\k),\Bbb C}{\tsize\bigwedge}^p(\k/\t)[p]$ is a K-projective
complex of weak $(\k,T)$-modules, since it is a translate of
$D(\ind_{\U(\k),\Bbb C}{\tsize\bigwedge}^p(\k/\t))$, and $D(P)$ is a
K-projective complex for any projective object $P$ in an abelian
category (see \cite{P2}, A1.4.5).  Now the functor $C_\t$ preserves
K-projectives since it is left adjoint to the forgetful functor from
equivariant complexes to weak complexes (see 2.5.2.(i)).  Therefore, by
3.2.4 we see that $Gr_pN(\k)$ is a K-projective
$(\k,T,N(\t))$-module. Hence we have proved

\proclaim{3.2.5. Corollary} If $T$ is reductive, then the graded objects
$Gr_pN(\k)$ corresponding to the Hochschild-Serre filtration are
K-projective $(\k,T,N(\t))$-modules.
\qed\endproclaim

We can now finish the proof of K-projectivity of $N(\k)$ and hence also of 3.1.5.

\proclaim{3.2.6. Theorem} If $T$ is reductive, $N(\k)$ is a K-projective 
$(\k,T,N(\t))$-module.
\endproclaim

\pf By 2.6.4, the only remaining thing to check is that for any $p$, 
$F_pN(\k)$ is isomorphic to $F_{p-1}N(\k)\oplus Gr_pN(\k)$ as a graded
$(\k,T,N(\t))$-module.

Let $\p$ be a $T$-invariant complement of $\t$ in $\k$; $\p$ exists since $T$ 
is reductive. Identifying $\k/\t$ with $\p$ gives an identification
$$
Gr_pN(\k)\cong (\U(\k)\otimes_\Bbb C{\tsize\bigwedge}^p(\p))[p]\otimes_\Bbb C 
{\tsize\bigwedge}(\t).
$$
As a graded module, the last module embeds into $F_pN(\k)$ via the inclusion 
induced by the isomorphism
$$
{\tsize\bigwedge}(\p)\otimes_\Bbb C{\tsize\bigwedge}(\t)@>>>{\tsize\bigwedge}(\k),
$$
which is given by multiplication.
 This inclusion clearly gives the required splitting.
\qed\enddemo

Now we will tie this with the well known construction of the relative standard 
complex used to calculate $(\k,T)$-cohomology. Note first that it follows from
2.2.6 that the counit morphism
$$
\eps:N(\k)@>>>\Bbb C,
$$
which is well-known to be a resolution of $\Bbb C$ by free $\U(\k)$-modules,
is also a K-projective resolution of the trivial $(\k,T,N(\t))$-module
$\Bbb C$. Namely, it is obvious that $\eps$ is a $T$-morphism and an 
$N(\t)$-morphism.

However, $\Bbb C$ is actually a $(\k,T)$-module, and we know from 2.5.2.(ii)
that the forgetful functor $Q$ from 
complexes of $(\k,T)$-modules into $(\k,T,N(\t))$-modules has a left adjoint 
$-\otimes_{N(\t)}\Bbb C$, that is, the functor of $N(\t)$-coinvariants. 
Therefore $\eps$ factors through a morphism
$$
\delta:N(\k)\otimes_{N(\t)}\Bbb C@>>>\Bbb C
$$
of complexes of $(\k,T)$-modules. Let us show that $\delta$ is a quasiisomorphism.
 Note that $\eps=\delta\circ p$, where
$p:N(\k)@>>>N(\k)\otimes_{N(\t)}\Bbb C$ is given by $n\mapsto n\otimes 1$; $p$
is actually the adjunction morphism. Therefore it is enough to show that $p$
is a quasiisomorphism. However, if we identify $N(\k)$ with 
$N(\k)\otimes_{N(\t)}N(\t)$, then $p$ gets identified with
$$
1\otimes\eps':N(\k)\otimes_{N(\t)}N(\t)@>>>N(\k)\otimes_{N(\t)}\Bbb C,
$$
where $\eps':N(\t)@>>>\Bbb C$ is the counit morphism, which is again a 
quasiisomorphism. By Poincar\'e-Birkhoff-Witt Theorem (2.2.1), $N(\k)$ is a free
right $N(\t)$-module for the right multiplication. Therefore, by \cite{BL1},
10.12.4.4, it is a K-flat $N(\t)$-module and it follows that $1\otimes\eps'$ is 
a quasiisomorphism.

On the other hand, since the functor $-\otimes_{N(\t)}\Bbb C$ preserves
K-projectives, $N(\k)\otimes_{N(\t)}\Bbb C$ is a K-projective complex
of $(\k,T)$-modules. In fact, its components are easily seen to be
$$
\U(\k)\otimes_{\U(\t)}{\tsize\bigwedge}^i(\p),
$$
so they are projective $(\k,T)$-modules. Namely, ${\tsize\bigwedge}^i(\p)$
are projective $T$-modules, that is, projective $(\t,T)$-modules, and
the above modules are $\ind_{\U(\k),\U(\t)}{\tsize\bigwedge}^i(\p)$, so they
are projective by 1.2.3. 

Note that the above modules are also 
finitely generated $\U(\k)$-modules.  To conclude 

\proclaim{3.2.7. Proposition} If $T$ is reductive, then the morphism 
$$
N(\k)\otimes_{N(\t)}\Bbb C @>>>\Bbb C
$$
constructed above is a resolution of the trivial $(\k,T)$-module $\Bbb C$
by projective $(\k,T)$-modules, which are finitely generated as 
$\U(\k)$-modules.
\qed\endproclaim
Of course, this is just the well-known relative standard complex (see 
e.g. \cite{BW}).

\subhead 3.3. Equivariant versus classical Zuckerman functors \endsubhead
In this section we relate the equivariant construction from Section 3.1 to 
the classical
Zuckerman construction. Let $V$ be an $(\A,T)$-module. Let $D(V)$ be the
corresponding complex of $(\A,T)$-modules (concentrated in degree 0). Consider 
the corresponding
$(\A,T,N(\t))$-module $Q(D(V))$ ($Q$ is as in 2.5.2.(ii)).
This is again a complex concentrated in
degree 0 and having $V$ as the $0^{\text th}$ component; the actions of $\A$
and $T$ are the same as on $V$, while the action of $N(\t)$ is trivial.
Let us consider the module $R(K,Q(D(V)))$. As a complex of weak $(\A,K)$-modules, 
with actions $\pi_R$ and $\nu_R$, it is just $D(R(K,V))$.
As a $(\k,T,N(\t))$-module, with action $\lambda$, it must be equal to 
$Q(D(R(K,V)))$ since it is concentrated in degree 0. It follows that
$$
\Gamma^{eq}_{K,T}(Q(D(V)))=\Homd_{(\k,T,N(\t))}(N(\k),Q(D(R(K,V)))).
$$
As we already mentioned, by 2.5.2.(ii), $Q:K(\M(\k,T))@>>>K(\k,T,N(\t))$ has a 
left adjoint, the functor $-\otimes_{N(\t)}\Bbb C$ of $N(\t)$-coinvariants. 
Since 
both these functors satisfy the conditions of 2.4.5, the stronger variant 
of adjunction from 2.4.5 holds, i.e., we have
$$
\Homd_{(\k,T,N(\t))}(N(\k),Q(D(R(K,V))))=\Homd_{(\k,T)}(N(\k)\otimes_{N(\t)}
\Bbb C,D(R(K,V))).
$$  
The actions of $\A$, $K$ and $N(\k)$ are given on the right hand side
in the same way as on the left hand side: $\A$ and $K$ act on the second 
variable, while $N(\k)$ acts on the first (by left multiplication). With 
this definition, it is obvious that the above equality is an equality of 
$(\A,K,N(\k))$-modules.

However, as we saw in 3.2.7, $N(\k)\otimes_{N(\t)}\Bbb C$ is a projective 
resolution of the trivial module $\Bbb C$ in the category of $(\k,T)$-modules.
It follows that as a complex of weak $(\A,K)$-modules,
$$
\Gamma^{eq}_{K,T}(Q(D(V)))=R\Hom_{(\k,T)}(D(\Bbb C),D(R(K,V))).
$$
In particular, we get

\proclaim{3.3.1. Lemma} Let $V$ be an $(\A,T)$-module. If $T$ is reductive, 
then the cohomology modules of $\Gamma^{eq}_{K,T}(Q(D(V)))$ are given by
$$
H^p(\Gamma^{eq}_{K,T}(Q(D(V))))=\Ext^p_{(\k,T)}(\Bbb C,R(K,V))
$$
as $(\A,K)$-modules. Here $R(K,V)$ is considered as a $(\k,T)$-module with 
respect to the $\lambda$-action.
\qed\endproclaim

Using 1.3.3, we finally get

\proclaim{3.3.2. Theorem} Assume that $T$ is reductive and that $\A$ is a flat
$\U(\k)$-module for the right multiplication (for example, this is true for 
$\A=\U(\g)$, where $(\g,K)$ is a Harish-Chandra pair). Let $V$ be an 
$(\A,T)$-module. Then the cohomology modules of $\Gamma^{eq}_{K,T}(Q(D(V)))$ 
are equal to the classical derived Zuckerman functors of $V$.
\qed\endproclaim

\enddocument
\end